\def\ZZ{{\bf Z}}

\outer\def\give#1. {\medbreak
             \noindent{\bf#1. }}                     
\outer\def\section #1\par{\bigbreak\centerline{\S
     {\bf#1}}\nobreak\smallskip\noindent}
\def\({\left(}
\def\){\right)}
\def\contract{{\,{\vrule height5pt width0.4pt depth 0pt}
{\vrule height0.4pt width6pt depth 0pt}\,}}

\def\sqr#1#2{{\vcenter{\hrule height.#2pt              
     \hbox{\vrule width.#2pt height#1pt\kern#1pt
     \vrule width.#2pt}
     \hrule height.#2pt}}}
\def\square{\mathchoice\sqr{5.5}4\sqr{5.0}4\sqr{4.8}3\sqr{4.8}3}
\def\qed{\hskip4pt plus1fill\ $\square$\par\medbreak}



\def\cA{{\cal A}}

\def\cC{{\cal C}}

\def\cE{{\cal E}}

\def\cG{{\cal G}}

\def\cL{{\cal L}}
\def\cM{{\cal M}}

\def\cO{{\cal O}}

\def\cR{{\cal R}}

\def\cT{{\cal T}}

 

\def\C{{\bf C}}
\def\cx#1{{\C}^{#1}}     
\def\cp1{{{\bf P}^1}}
\def\R{{\bf R}}

\def\Z{{\bf Z}}
\def\bar{\overline}              
 
 
\def\bd{\partial}                

\magnification=\magstep1

\centerline{Polynomial Diffeomorphisms of $\cx 2$:}
\centerline{VI. Connectivity of $J$}
\bigskip
\centerline{Eric Bedford and John Smillie}
\bigskip
\section 0. Introduction

Polynomial maps $g:\cx{}\to\cx{}$ are the simplest holomorphic maps with
interesting dynamical behavior. The study of such maps has had an important influence on
the field of dynamical systems.  On the other hand the traditional focus of the field of
dynamical systems has been in a different direction: invertible maps or diffeomorphisms. Thus
we are led to study polynomial diffeomorphisms
$f:\cx2\to\cx2$, which are the simplest holomorphic diffeomorphisms
with interesting dynamical behavior.   

Two features are apparent in much of the contemporary work on polynomial maps
of $\C$ (cf.\  [DH]).  The first is a focus on the connectivity of the Julia
set. The second is the use of computer pictures as a guide to research.
Computer pictures do not substitute for proofs but they have  provided a  tool
that has been used to guide research. In this paper we consider these ideas in
the context of polynomial diffeomorphisms of
$\C^2$. This approach to the study of polynomial diffeomorphisms originates with Hubbard.

 For a polynomial map of $\C$ the ``filled Julia set'' $K\subset\cx{}$ is the set of points
with bounded orbits, and the Julia set $J$ is defined to be the boundary of $K$.
The Julia set has several analogs for diffeomorphisms of $\C^2$. Since
$f:\cx2\to\cx2$ is invertible, we can distinguish properties of points
based on both forward and backward iteration.   The sets $K^+$ (resp.\
$K^-$) consist of points with bounded forward (resp.\ backward) orbits under
$f$.  We write $U^\pm$ for the complementary sets $U^\pm=\C^2-K^\pm$. The set of
points whose orbits are bounded in both forward and backward time is $K=K^+\cap
K^-$. The sets $J^\pm:=\partial K^\pm$ are analogues of the Julia set, as is the
set $J=J^+\cap J^-$. We use the notation $J^-_+$ for $J^-\cap U^+$. We will see that in
some cases the set  $J^-_+$ plays the role of the Fatou set for polynomial maps of $\C$.
The focus of this paper is to investigate  the 
$J$-connected/ $J$-disconnected dichotomy in the case of polynomial diffeomorphisms of
$\cx2$ and relate it to the structure of the sets $J^-\cap U^+$ and $J^+\cap U^-$.

  One of the attractive features of the study of polynomial maps of $\C$ is that the Julia
sets can be drawn by computer. Thus the connectivity properties of the Julia set can often be
demonstrated (visually if not rigorously) by means of computer pictures.  One of the daunting
features of the study of polynomial diffeomorphisms of $\C^2$ is that the sets of fundamental
importance are complicated subsets of $\C^2$.  As a consequence of our investigations we will
show that it is possible to ``see'' the connectivity of the Julia set $J$ for polynomial
diffeomorphisms of $\C^2$. 

 Hubbard has suggested the following computer experiment. Let $f$ be a polynomial diffeomorphisms of
$\cx2$ and $p$ be a
 periodic saddle point of $f$.  The unstable manifold of $p$,
$W^u(p)$, is an immersed submanifold conformally equivalent to $\C$. We have a
partition of $W^u(p)$ into subsets $W^u(p)\cap K^+$ and $W^u(p)\cap U^+$. (We can
view this as a partition of $\C$.)  This partition is easily drawn by
computer, and  the sets $W^u(p)\cap K^+$ have many local features of  Julia sets. 

While the unstable manifolds are easily drawn and yield complicated pictures, these pictures
 depend on choice of the point $p$.  In order to make these computer pictures a
useful tool for investigating the dynamics of $f$ it is important to know which features
of these pictures reflect dynamically significant properties of the diffeomorphism.  In
particular it is useful to know which features are independent of the choice of periodic
point. We say that
$f$ is {\sl unstably connected with respect to the point $p$ } if some component of
$W^u(p)\cap U^+$ is simply connected. The following Theorem shows that this property is
independent of the point $p$ and furthermore that it is equivalent to the existence of a
certain geometric structure of $J^-_+$.

\proclaim Theorem 0.1.  The following are equivalent:
\item{1.}  For some periodic saddle point $p$, some component of $W^u(p)\cap
U^+$ is simply connected.
\item{2.}  The set $J^-_+$ has a lamination by simply connected leaves so
that for any periodic saddle point $p$ each component of $W^u(p)\cap U^+$ is a
leaf of this lamination.
\item{3.}  For any periodic saddle point $p$, each component of $W^u(p)\cap
U^+$ is simply connected.

This result is
a consequence of Theorems 2.1 and 4.1. 

If the equivalent conditions of the Theorem hold, we say that $f$ is
{\sl unstably connected}.
When $f$ is a
polynomial diffeomorphism, the inverse function is also a polynomial
diffeomorphism. Replacing $f$ by $f^{-1}$ interchanges the sets $J^-$ and
$J^+$, so any result for, or property of, $J^-$ has an analog for $J^+$.  We
say that $f$ is {\sl stably connected} if the equivalent properties of the
previous theorem hold for
$f^{-1}$. In \S5 we show that the connectivity of $J$ is determined by the  stable/unstable
connectivity of $f$.

\proclaim Theorem 5.1.  The set $J_f$ is connected if and only if the diffeomorphism
$f$ is either stably connected or unstably connected. 

 It follows from Theorem 5.1 that
the connectivity of the Julia set $J$ of a polynomial diffeomorphism of
$\C^2$ can be determined ``empirically''.  It suffices to pick a periodic saddle point
and observe the escape locus inside the stable and unstable manifolds, and according to
Theorem 0.1 and its analog for stable connectivity the topology of these sets determines
the connectivity of $J$.

In one variable, questions involving the connectivity of $J$ are bound up with
properties of the critical points of $g$.  In particular $J$ is connected if
and only if there are no critical points with unbounded orbits.  Now a polynomial
diffeomorphism by definition has no critical points, and it seems difficult to define a single
 analog of critical points which works in all settings.  On the other hand, an analog of
critical points with unbounded forward orbits was described in [BS5]. This 
is the set
$\cC^{u}$ of critical points of the Green function $G^+$ restricted to unstable manifolds. 

\proclaim Theorem 7.3.  $f$ is unstably connected if and only if for $\mu$ almost
every point $p$, $W^u(p)\cap U^+$ contains no unstable critical points.

Recall that the Jacobian determinant of a polynomial diffeomorphism $det D(f)$ is
constant, and $det D(f^{-1})=det D(f)^{-1}$.  Replacing $f$ by
$f^{-1}$ if necessary, we may assume that $|det Df|\le1$.  We say
that $f$ is dissipative if $|det Df|<1$ and volume preserving if $|det Df|=1$. 

Combining the previous theorem with results from [BS5] we have:

\proclaim Corollary 7.4. If $f$ is dissipative then $f$ is not stably connected. If $f$ is
volume preserving, then $f$ is stably connected if and only if $f$ is unstably connected.

Combined with Theorem 5.1 this immediately yields:

\proclaim Theorem 0.2. Assume $|det Df|\le 1$.  Then $J$ is connected if and only if $f$ is
unstably connected.

Earlier we described a method of determining the connectivity of
$J$ experimentally by considering the stable and unstable manifolds of a periodic
point.  It follows from Theorem 0.2 that we can determine the connectivity of $J$ with only
half as much data.

When $f$ is unstably connected we can give further information about the lamination of
$J^-_+$.  Our result can be compared with the fact that for polynomial maps in one variable
the connectivity of the Julia set is equivalent to the existence of a canonical model for the
dynamics on the complement of the Julia set. 

The complex solenoid is defined as the set of bi-infinite sequences
$\Sigma_+=\{z=(z_j):|z_j|>1,j\in\ZZ,z^d_j=z_{j+1}\}$ with the product topology.  The
induced mapping
$\sigma:\Sigma_+\to\Sigma_+$, obtained by shifting sequences to the left, is a
homeomorphism.  Furthermore,
$\Sigma_+$ has a natural lamination by Riemann surfaces.

\proclaim Theorem 3.2.  If $f$ is unstably connected, then there is a continuous
map $\Phi:J^-_+\to\Sigma_+$ which semi-conjugates $f|_{J^-\cap U^+}$ to the
mapping $\sigma|_{\Sigma_+}$, i.e.\ $\sigma\circ\Phi=\Phi\circ f$. 
Furthermore, $\Phi$ preserves the lamination structure and is a holomorphic
bijection on each leaf.

Given this laminar structure on $J^-_+$  we
 define external rays  to be gradient lines of the function $G^+$ restricted to
leaves of the lamination. As in the case of polynomial maps of $\C$ it is interesting
to know when these rays ``land,'' i.e., converge to points in $J$.
 Even in one variable it is too much to ask that all rays land, unless we impose
additonal hypotheses on the map. It is known that the set of rays which fail to land has
Lebesgue measure zero. Theorem 3.1 shows that for polynomial diffeomorphisms, there is a natural
measure on the set of rays, with respect to which almost all rays land. In [BS7] we show that
when
$f$ is unstably connected and hyperbolic all rays land.  We also show in [BS7] that
with the same hypotheses the semiconjugacy to the solenoid can be replaced by a
conjugacy.

We end this section with a guide to the organization of the paper.
As a tool for relating unstable connectivity at a point and the laminar
structure on $J^-_+$, we introduce in \S2 a condition
(\dag), and we show that the existence of a laminar structure is a consequence of (\dag). In
section
\S3 we show that (\dag) implies the existence of a semiconjugacy to the
solenoid. In \S4 we show that the topological condition of unstable connectivity implies
(\dag). In \S5 we relate stable and unstable connectivity to the connectedness of
$J$. In \S6 we collect several properties equivalent to unstable connectivity. These are
summarized in Theorem 6.9. In particular this theorem shows that unstable connectivity
can be characterized without making reference to unstable manifolds.  In
\S7 we establish a relationship between unstable connectivity and critical points. 

\section 1.  Notation and Preliminaries

In this section we briefly review standard terminology and facts about polynomial
diffeomorphisms of $\C^2$. For a discussion at greater length see \S1 of
[BS5].

 We consider mappings of the form
$f=f_1\circ\cdots\circ f_m$, where each
$f_j$ has the form
$$f_j(x,y)=(y,p_j(y)-a_jx),$$
and $p_j(y)$ is a monic polynomial of degree $d_j$, and $a_j\in\cx{}$ is
nonzero.  It follows that $f$ has degree $d=d_1\cdots
d_m$, and the $n$-fold iterate $f^n=f\circ\cdots\circ f$ has degree
$d^n$.  Let $K^\pm$ denote the points which remain bounded in
forward/backward time, and  set $K=K^+\cap K^-$, $J^\pm=\partial K^\pm$,
and $J=J^+\cap J^-$.   Define the projections
$\pi_x(x,y)=x$ and $\pi_y(x,y)=y$, and it follows that
$\pi_x(f^n)=y^{d^{n-1}}+\cdots$ and 
$\pi_y(f^n)=y^{d^n}+\cdots\,$.  Set    
$$G^\pm(x,y)=\lim_{n\to+\infty}{1\over d^n}\log^+\left\Vert f^{\pm
n}(x,y)\right\Vert =\lim_{n\to+\infty}{1\over d^n}\log^+\left|\pi_y\circ f^{\pm
n}(x,y)\right|,$$ 
then $G^\pm$ is continuous and plurisubharmonic on $\cx2$,
$K^\pm=\{G^\pm=0\}$,  and $G^\pm$ is pluriharmonic on the sets
$U^\pm:=\cx2-K^\pm$.  The currents $\mu^+$ and $\mu^-$ are defined by
$\mu^\pm=(1/2\pi)dd^cG^\pm$. The harmonic measure is given by $\mu=\mu^+\wedge\mu^-$.

For $R>0$  define
$V^+=\{|y|>|x|,|y|>R\}$.  It follows that one may choose
$R_0$ sufficiently large that
$$J^-\cap\{G^+>R_0\}\subset V^+.\eqno(1.1)$$
If $R$ is chosen sufficiently large, then for all $n\ge0$, one may choose a
$d^n$-th root of
$\pi_y\circ f^n$ on $V^+$ such that $(\pi_y(f^n))^{1\over d^n}\approx y,$ and
with this choice of root, define
$$\varphi^+(x,y)=\lim_{n\to+\infty}(\pi_y\circ f^n(x,y))^{1\over d^n},$$
which is analytic on $V^+$ and satisfies $\varphi^+\circ
f=(\varphi^+)^d$ (see [HOV]).  Further, 
$$\varphi^+(x,y)=y+O(1),\hbox{\ and\ \ }{\partial\varphi^+(x,y)\over \partial
y}=1+O(|y|^{-1})\eqno(1.2)$$ on $V^+$.  The uniformity of the $O$ terms means
that for fixed $x$,
$y\mapsto\varphi^+(x,y)$ is univalent for $|y|\ge\max(R,|x|)$.  Thus, if we set
$s=\varphi^+$, then $(x,s)$ is a global coordinate system on $V^+$.

Let $\nu$ be an ergodic invariant measure supported on $J$. Associated to $\nu$ are two
Lyapunov exponents $\lambda^+(\nu)$ and $\lambda^-(\nu)$ which describe the growth of tangent
vectors under $Df^n$. We say that $\nu$ is an index one hyperbolic measure (or simply a
hyperbolic measure when no confusion will result) if
$\lambda^-(\nu)<0<\lambda^+(\nu)$. If $\nu$ is an index one
hyperbolic measure then the following conditions hold:
\item{(1)} for
$\nu$ a.e.\ $p$ there are complex one dimensional linear subspaces  $E^s_p$ and $E^u_p$ of the
tangent space $T_p\cx2$ such that $E^s_p\oplus E^u_p=T_p\cx2$
where
 the families $\{E^{s/u}_p\}$ are
invariant: $Df_p(E^{s/u}_p)=E^{s/u}_{fp}$. 
\item{(2)} for $\nu$
a.e.\ point $p$, we have 
$$\lim_{|n|\to\infty}{1\over n}\log\Vert
Df^n_p|_{E^{s/u}_p}\Vert=\lambda^{-/+}(\nu).\eqno(1.3)$$ 

Conversely if there are numbers $\lambda^-<0<\lambda^+$ such that (1) and (2) hold then $\nu$
is an index one hyperbolic measure with exponents $\lambda^\pm$.
 The concept of hyperbolic measure will be used frequently
throughout this paper; see  [P] or [KM] for further information.

Examples of index one hyperbolic measures include the  measure given by
the average of point masses over the orbit of a periodic saddle point. Such periodic saddle
points exist in profusion for every
$f$ (see [BLS2]). Another important example of a hyperbolic measure is the harmonic measure
$\mu$ (see [BS3]).

A consequence of the Pesin theory is that there is a set $\cR\subset J$ which has full
$\nu$ measure for any index one hyperbolic measure $\nu$ and such that for $p\in\cR$ the
stable/unstable manifolds
$$W^{s/u}(p) =\{q\in\cx2:\lim_{n\to\pm\infty}\hbox{\rm dist}
(f^nq,f^np) =0\}$$
are smooth submanifolds of $\cx2$.  In fact, (see [BLS], [W]) these are Riemann
surfaces which are conformally equivalent to $\cx{}$.  In general these
submanifolds vary only measurably with $p$. For our purposes we may {\it define} $\cR$ to be
the set of points $p\in J$ for which $W^{s/u}(p)$ are Riemann
surfaces conformally equivalent to $\cx{}$.

\section 2.  Extension of $\varphi^+$ and Laminar Properties of $J^-_+$

In this section we describe a technical condition on complex disks contained
in $J^-_+$ which will play an important role in this paper. We will show that when
there is a disk that satisfies this condition, 
$J^-_+$ has a lamination by complex leaves. In \S4 we will give
topological hypotheses that imply the existence of such a disk.

 Let $M$ be a Riemann surface.  A positive harmonic
function on
$M$ will be called {\it minimal} if it generates an extreme ray in the cone of
positive harmonic functions.  In other words, if a minimal harmonic function $h$
can be written as a sum of positive harmonic functions $h=h_1+h_2$ on $M$, then
$h,h_1$, and $h_2$ are all constant multiples of each other.  In case
$M=\Delta$ is the unit disk, the minimal harmonic functions are just the
positive constant multiples of the Poisson kernel
$$P(z,e^{i\theta})=\Re\({e^{i\theta}+z\over e^{i\theta}-z}\)$$
for some real $\theta$.
We note that under composition with an automorphism of the disk $\Delta$, the
Poisson kernel $P(z,1)$ is taken to $cP(z,e^{i\theta})$ for some real
numbers $c>0$ and $\theta$.  For any positive harmonic function $h$ on the unit
disk, there is a unique positive Borel measure $\lambda$ on $\partial\Delta$
for which we have the Herglotz representation
$$h(z)=\int P(z,e^{i\theta})\lambda(\theta).$$

It turns out that minimality of $G^+$ will enter frequently in our work in a
rather specific context.  So we formalize this.  We call $\cO\subset\cx2$ a {\it
complex disk} if there is a holomorphic injective immersion  of the
unit disk
$\psi:\Delta\to\cx2$ with
$\cO=\psi(\Delta)$.  We say that a complex disk $\cO$ satisfies condition (\dag) if the
following three properties hold:

$$\eqalign{&\cO\subset J^-_+, \quad G^+|\cO  \hbox{\rm\ is\  minimal,\ and\  }\cr
\hbox{\rm\  for\ each\ }&j\in\Z \hbox{\rm\   either \ }
\cO\cap f^j\cO=\emptyset  \hbox{\rm\  or\ }\cO=f^j\cO.\cr}\eqno(\dag)$$
\smallskip

In the most typical situation $\cO$ will be a component of $W^u(p)-K^+$ where
$W^u(p)$ is an unstable manifold of a point (not necessarily periodic). In this
case the the minimality of $G^+$ is the only part of (\dag) that is not automatic.
We say that a polynomial diffeomorphism satisfies condition (\dag) if it posesses a disk which
satisfies condition (\dag).

Let us recall the definition of a Riemann surface lamination of a topological
space $X$ (cf. [C]).  A chart is a choice of an open set $U_j\subset X$, a
topological space $Y_j$, and a map
$\rho_j:U_j\to\C\times Y_j$ which is a homeomorphism onto its image.  An atlas is
a collection of charts such that $\{U_j\}$ covers $X$.  The set of points of $U_j$
for which the second coordinate of
$\rho_j$ assumes a fixed value is called a plaque.  For coordinate charts
$(\rho_i,U_i,Y_i)$ and $(\rho_j,U_j,Y_j)$ with $U_i\cap U_j\ne\emptyset$, the
transition function is the homeomorphism from $\rho_j(U_i\cap U_j)$ to
$\rho_i(U_i\cap U_j)$ defined by $\rho_{ij}=\rho_i\circ\rho_j^{-1}$.  A Riemann
surface lamination of a topological space $X$ is determined by an atlas of charts
which satisfy the following consistency condition: the transition functions may be
written in the form $\rho_{ij}=(g(z,y),h(y))$, where for fixed $y\in Y_j$ the
function $z\mapsto g(z,y)$ is holomorphic.  The condition on the transition
functions gives a consistency between the plaques defined in $U_j$ and those in
$U_i$.  Thus plaques fit together to make global manifolds called leaves of the
lamination, and each leaf has the structure of a Riemann surface.

One of the equivalent conditions in Theorem 2.1 below is that $J^-_+$ carries a
unique Riemann surface lamination.  We  remark after Proposition 2.7 that this
lamination in fact carries a special affine structure.

The following theorem collects the basic results of this section.

\proclaim Theorem 2.1.  If there is a complex disk $\cO$ satisfying (\dag), then
\item{(1)}  $\varphi^+$ extends to a continuous mapping
$\varphi^+:J^-_+\to\{|\zeta|>1\}$   which satisfies the functional equation 
$$\varphi^+(f(p))=(\varphi^+(p))^d.$$
\item{(2)} $J^-_+$ has a Riemann surface lamination $\cM^-$ with $\cO$ as one leaf.
\item{(3)}  For each leaf $M$ of $\cM^-$, $G^+|M$ has no critical points.  
\item{(4)}   For each leaf $M$ of the lamination $\cM^-$, the restriction
$\varphi^+|_M:M\to\{|\zeta|>1\}$ is a holomorphic covering map.
\item{(5)}  Each leaf $M$ of $\cM^-$ is a disk satisfying (\dag).
\item{(6)}  Each leaf of $\cM^-$ is dense in $J^-_+$.
\item{(7)}  $\cM^-$ is the unique lamination of $J^-_+$ by Riemann surfaces.

The remainder of this section is devoted to the proof of Theorem 2.1. In particular we will
make the standing assumption in this section of the paper that there exists a disk $\cO$ which
satisfies (\dag).  The general plan of the proof is as follows.  Let $\cO^\infty$ be defined as
$\bigcup_{j\in\Z}f^j\cO$ in $J^-_+$. Let $\hat\cO$ denote the closure in $U^+$ of the set
$\cO^\infty$.  We begin by proving analogs of (1--5) with $J^-_+$ replaced
by $\hat\cO$.
In Theorem 2.6 we prove an analog of (1).  In Proposition 2.7 we prove an analog of
(2). In Corollary 2.8 we prove an analog of (3). In Corollary 2.9 we prove an analog of
(4). In Corollary 2.10 we prove (5).  Corollary 2.14  proves (6) and  
shows that $\hat\cO=J^-_+$. Thus it follows that items (1--5) hold for $J^-_+$ as
claimed.  (7) is proved in Corollary 2.18.

For any set  $X\subset J^-_+$ we will use the notation:
$$X(\rho):=X\cap\{G^+\ge\log\rho\}.$$
We saw in \S1 that $\varphi^+$ is defined for points in $J^-_+$ with $G^+$
suffficiently large. Let $\rho_0$ be a fixed constant chosen large enough that
$\varphi^+$ is defined on  $J^-_+(\rho_0)$.

We will also use the notation $\C(\rho):=\{\zeta\in\C:|\zeta|\ge\rho\}$. With this
notation it follows that for
$\rho\ge\rho_0$ the map $\varphi^+$ maps the set $J^-_+(\rho)$ to the set
$\C(\rho)$.

\proclaim Proposition 2.2.  Let $H=\{x+iy:x>0\}$ be the right half-plane, and let
$\cO$ be a disk satisfying (\dag).  We can choose a conformal coordinate
$\alpha:H\to\cO$ such that for $x>\log \rho_0$ we have $\varphi^+(\alpha(x+iy))$ is
defined and 
$\varphi^+(\alpha(z))=e^z$.

\give Proof.  By hypothesis, there is an imbedding
$\psi:\Delta\to\cO\subset\cx2$ such that $g=G^+\circ\psi$ is a positive
multiple of the Poisson kernel function with pole at
$e^{i\kappa}\in\partial\Delta$.  Choose a conformal map from the right half
plane $H$ to $\Delta$ which sends the point at infinity to $e^{i\kappa}$.  Let
$\alpha:H\to\cO$ be the composition of map from $H$ to $\Delta$ and $\psi$.  On
the right half-plane the Poisson kernel function with pole at infinity is
$x+iy\mapsto x$.  Thus $G^+(\alpha(x+iy))=cx$. By composing $\alpha$ with
multiplication by a scalar we may assume that $c=1$.  If
$x>\log\rho_0$ then $\varphi^+(\alpha(x+iy))$ is defined, and
$G^+(\alpha(x+iy))=\log|\varphi^+(x+iy))|$.  Since $\log|e^{x+iy}|=x$, we have
$\log|\varphi^+(x+iy)|=\log|e^{x+iy}|$ so that
$$\log|\varphi^+(x+iy)/e^{x+iy}|=0$$
$$\Re(\varphi^+(x+iy)/e^{x+iy})=1.$$
A holomorphic function with constant real part is constant so  that
$\varphi^+(x+iy)=e^{i\theta_0}e^{x+iy}$.  Now by composing
$\alpha$ with the translation $z\mapsto z-i\theta_0$ we have the required
parametrization. \qed

  Since  the map
$exp:\{x+iy: x>\log\rho\}\to\{\zeta:|\zeta|>\rho\}$ is a covering when $\rho\ge\rho_0$,
we have:

\proclaim Corollary 2.3. If $\cO$ satisfies $(\dag)$ then
$\varphi^+|_{\cO(\rho)}:\cO(\rho)\to\{|\zeta|>\rho\}$ is a covering
map   when $\rho>\rho_0$.

For any subsets $X\subset J^-_+$ and  $E\subset\{G^+>\rho\}$ with $\rho\ge\rho_0$,
we will use the notation
$X_E:=X\cap(\varphi^+)^{-1}E$.

\proclaim Lemma 2.4.  If $G\subset\{|s|>\rho_0\}$ is simply connected and $\zeta\in
G$, then there  a homeomorphism
$H: G\times\hat\cO_{\{\zeta\}}\to \hat\cO_G$ satisfying 
$\varphi^+(H(s,p))=s$. Furthermore $H(G\times\cO^\infty_{\{\zeta\}})=\cO^\infty_G$.

\give Proof.  Condition (\dag) implies that $\cO^\infty$ can be written as a disjoint
union of disks of the form 
$f^n(\cO)$, which again satisfy (\dag).  For
$p\in\cO^\infty$ let $\cO'$ denote the (unique) disk containing $p$, and let $L_p$ be
the component of
$\cO'_G$ which contains $p$.
 Now by Corollary 2.3 $\varphi^+|_{L_p}$ is a covering map with connected
domain and simply connected range, so in fact $\varphi^+|_{L_p}$ is a
bijection.  Let
$g_p:G\to L_p\subset\cx2$ be the inverse map, so $g_p$ is a holomorphic.  By
the remark following (1.2), the functions $(\pi_x,\varphi^+)$ form a coordinate
system on $V^+$, so there is a function $h$ such that $g_p(z)=(h(z,p),z)$.  Further, if
we set
$E:=\cO^\infty_{\{\zeta\}}$, then in these coordinates we may identify $E$ with a
subset of $\C$.

 Now the function $h:G\times E\to\cx{}$ is a
holomorphic motion, which is to say that: 
\item{(1)} $z\mapsto h(z,p)$ is holomorphic, and
\item{(2)} $p\mapsto h(z,p)$ is injective.  

The second property follows from the fact that
$\cO^\infty$ consists of subsets of pairwise disjoint disks in $\cx2$ which are given
as graphs $\zeta\mapsto(h(\zeta,p),\zeta)$ in the $(\pi_x,\varphi^+)$-coordinates.  
There is a well-developed theory of holomorphic motions which shows that
holomorphic motions can be extended to quasiconformal maps
$\cx{}\ni t\mapsto h(z,t)$.  We need only the basic observation from [MSS] that $h$
has an extension to a function $\hat h: G\times\bar {\cO^\infty_{\{\zeta\}}}\to\cx{}$,
which is also a holomorphic motion, and which is continuous.  Now
$\bar{\cO^\infty_{\{\zeta\}}}= {\hat\cO}_{\{\zeta\}}$. We define
$H:G\times\hat\cO_{\{\zeta\}}\to\hat X_G$ by the formula
$H(t,z)=(t,\hat h(t,z))$, where the right hand side is given in the
$(\pi_x,\varphi^+)$ coordinates on $\hat X_G$.  The function
$H$ is a continuous surjection.  The fact that
$\hat h$ satisfies (2) shows that $H$ is injective.   $H$ is proper, so it is a
homeomorphism on each compact set. This implies that $H$ is a homeomorphism.
\qed

Next we show that the function $\varphi^+$ has a canonical extension to all of
$\hat\cO$. This proof will use a purely topological fact which we prove first.

Let $\pi:\tilde Y\to Y$ be a covering space. Let $\psi: X\to Y$ be a map. A {
lift} of
$\psi$ is a map
$\psi': X\to\tilde Y$ such that $\pi\circ \psi'=\psi$. We will reduce the problem of
extending
$\varphi^+$  to a problem of finding a lift. On well-behaved spaces lifting
problems can be translated into problems in terms of the fundamental group. The
spaces we are dealing with here are not locally connected so this translation is not
possible. On the other hand standard techniques of topology can be used to solve the
problem.  Let $A$ be a closed subset of the topological space $X$.
 To say that $A$ is a strong deformation retract of $X$ means that we are given a
function $F:X\times I\to X$ (the retraction function) such that:
\item{(1)} $F(x,0)=x$ for $x\in X$
\item{(2)} $F(x,1)\in A$ for $x\in X$
\item{(3)} $F(x,t)=x$ for $x\in A$.

\proclaim Lemma 2.5. Let $A\subset X$ be a strong deformation retract. 
Let $\pi:\tilde Y\to Y$ be a covering map. Let $\psi: X\to Y$
be a map and assume we are given a map $\psi_1: A\to \tilde Y$ which is a lift of
$\psi|A$. Then there is a unique lift $\psi'$ of $\psi$ which agrees with $\psi_1$ on
$A$.

\give Proof. Let $F$ denote the strong retraction function as above, and define 
$G:X\times I\to X$ by  $G(x,t)=\psi\circ F(x,1-t)$. The function
$x\mapsto
\psi_1\circ F(x,1)$ is well defined since $F(x,1)\in A$ and it is a lift of the function
$x\mapsto G(x,0)$. The homotopy lifting property of the covering map $\pi$
[S, p. 67, Theorem 3] gives us a unique lift $G'$ of $G$ with the property that
$G'(x,0)=\psi_1\circ F(x,1)$. The restriction of $G'$ to the set $A$ is a lift of a constant
homotopy (independent of $t$). By uniqueness of lifts of paths [S, p. 68,
Theorem 5] the restriction of
$G'$ to
$A$ is itself a constant homotopy. It follows that $G'(x,t)=\psi_1(x)$ for $x\in A$. Now
if we set $\psi'(x)=G(x,1)$ then $\psi'$ is a lift of $\psi\circ G(x,1)=\psi$ and
$\psi'(x)=x$ for $x\in A$.

Now if $\psi''$ is any other lift of $\psi$ which agrees with $\psi_1$ on $A$ then
$G''=\psi''\circ F$ is a lift of $G$ with the property that $G''(x,0)=\psi_1\circ F(x,1)$. By
the uniqueness property of lifts of homotopies we have $G''=G'$ hence
$\psi''(x)=G''(x,1)=G'(x,1)=\psi'$.   \qed

\proclaim Theorem 2.6. If $f$ satisfies $(\dag)$ then the function $\varphi^+$ has a continuous
extension to all of
$\hat\cO$ which satisfies the functional equation 
$$\varphi^+(f(p))=(\varphi^+(p))^d.\eqno(2.1)$$
There is a unique  extension satisfying (2.1).

\give Proof. The function $\varphi^+$ is defined on $\hat\cO(\rho_0)$ and satisfies
(2.1) for $p\in \hat\cO(\rho_0)$. To prove the Proposition, it suffices to show that
$\varphi^+$ has a unique, continuous extension to $\hat\cO(\rho_0^{1\over d})$,
which satisfies (2.1).  Repeating this argument allows us to extend the function
successively to sets $\hat\cO(\rho_0^{1/d^n})$ and thus to $\hat\cO$ which is the
union of these sets.  

Let us define
$\pi:\C \to\C$ by
$\pi(z)=z^d$  and write
$\psi=\varphi^+\circ f$. Finding a $\psi'$ which satisfies (2.1) is
equivalent to finding  $\psi':\hat\cO(\rho^{1/d})\to\C-\bar\Delta$ such that
$\pi\circ
\psi'=\psi$.
 The function  $\pi:\C-\bar\Delta
\to\C-\bar\Delta$ is a covering map. This is a problem of finding a lift $\psi'$ with the
added constraint that we want $\psi'$ to have prescribed values on $\hat\cO(\rho_0)$. 

In order to apply  Lemma 2.5, let $A=\hat\cO(\rho)$, let $X=\hat\cO(\rho^{1/d})$, and
let $\psi_1=\varphi|_{\hat\cO(\rho)}$. To verify the hypotheses of the  Lemma it
suffices to show the following. For
$1<\alpha<\beta$ the set
$\hat\cO(\beta)$ is a deformation retract of
$\hat\cO(\alpha)$.

   Since $f^n(\hat\cO(\rho))=\hat\cO(\rho^{d^n})$ by applying a sufficiently high
power of $f$ we may assume that $\rho_0\le\alpha<\beta$. Applying the
homeomorphism $f^n$ does not change the topological properties of the sets. We we
will construct a deformation retraction by using
$\varphi^+$ to lift a deformation retraction of
$\{|z|\ge\alpha\}$ to $\{|z|\ge\beta\}$. 

 Let $F:\{|z|\ge\alpha\}\times I\to \{|z|\ge\beta\}$ be a strong deformation
retraction.  Thus $F(z,0)=z$, $|F(z,1)|\ge\beta$ and $F(z,t)=z$ for $|z|\ge\beta$.
We may assume that
$F$ preserves radial lines. Define $G: \hat\cO(\alpha)\times I\to \C(\alpha)$ by
$G(p,t)=F(\varphi^+(p),t)$. We wish to find a function $G': \hat\cO(\alpha)\times I\to
\hat\cO(\alpha)$ which satisfies $\varphi^+\circ G'=G$ and
 $G'(p,0)=p$. Since $\varphi^+$ is a covering on each leaf the homotopy lifting
property of covering spaces gives us such a unique such function. A priori all we
know is that this function is continuous when restricted to each leaf. To see that
$G'$ is continuous on $\hat\cO(\alpha)$ we use the product structure given by Lemma
2.4. Consider a set $S=\{|z|\ge\alpha,\ \theta_0\le\arg(z)\le\theta_1\}$. Since $S$ is
simply connected the set $\hat\cO_S$ has a product representation as
$S\times\hat\cO_{\{\zeta\}}$ by Lemma 2.4. Now the restriction of $F$ to $S$ is a
deformation retraction as is the restriction of $G$ to $\hat\cO_S$.  To prove that $G$
is continuous it suffices to prove that $G$ restricted to $\hat\cO_S$ is continuous for
any such set
$S$. Now we can use the product structure to define a deformation retraction on
$\hat\cO_S$ by the formula $(z,w)\mapsto(F_t(z),w)$ and it is a lift of $F$. This 
function is clearly continuous and by the uniqueness of lifts it must be equal to $G$
restricted to
$\hat\cO_S$.  \qed

We now define the lamination $\cM^-$ on the set $\hat\cO$.

\proclaim Proposition 2.7. There is a Riemann surface lamination of $\hat\cO$. The
charts have the form $\rho_G:\hat\cO_G\to G\times \hat\cO_{\{\zeta\}}$ and satisfy
the condition $\pi_1\circ\rho_G=\varphi^+$.

\give Proof. For each point $p\in\hat\cO$ we will construct a set $G$ and a chart as
above. If $G^+(p)>\log\rho_0$ then let $G\subset \{|z|>\rho_0\}$ be a simply
connected open set that contains $\varphi^+(p)$. Let $\rho_G$ be the inverse of the
function $h$ given by Lemma 2.4.

For a  $p\in\hat\cO$ not of the above form there is an $n>0$ so that
$|\varphi^+(f^n(p))| >\rho_0$. Let $\xi=\varphi^+(f^n(p))$. Let
$\zeta_1,\ldots,\zeta_{d^n}$ be the roots of $z^{d^n}=\xi$ so that
$\zeta_1=\varphi^+(p)$. Let $H$ be a simply connected subset of $\{|z|>\rho_0\}$ that
contains $\xi$. Let $G_1,\ldots, G_{d^n}$ be the components of $\{z: z^{d^n}\in
H\}$ numbered so that $\zeta_k\in G_k$. Now
$f^{-n}(\hat\cO_H)=\hat\cO_{G_1}\cup\ldots \cup\hat\cO_{G_{d^n}}$ where the sets
on the right-hand side are disjoint. We have a map $\rho_G:\hat\cO_H\to H\times
\hat\cO_{\{\xi\}}$.  We can write $\hat\cO_{\{\xi\}}$ as
a disjoint union $\hat\cO_{\{\zeta_1\}}\cup\ldots\cup
\hat\cO_{\{\zeta_{d^n}\}}$. In fact
$f^{-n}\rho_G^{-1}(H\times\hat\cO_{\{\zeta_k\}})=
\hat\cO_{G_k}$. This is because for each plaque $P$ the function $\varphi^+\circ
f^{-n}\circ(\varphi^+|_P)^{-1}$ is a branch of the $d^n$-th root function so its
values lie in $G_k$ if and only if its value at $\xi$ is $\zeta_k$. So the function
$\rho_H\circ f^n:\cO_{G_1}\to H\times
\hat\cO_{G_1}$ gives a coordinate chart on $\cO_{G_1}$.

It remains to check the form of the overlap functions. If $\hat\cO_G'$ and
$\hat\cO_G''$ intersect, then their intersection is $\hat\cO_H$ where $H=G'\cap
G''$. We can analyze the transition function in terms of the functions from
$\hat\cO_H$ to $\hat\cO_G'$ and from
$\hat\cO_H$ to $\hat\cO_G''$.  Thus we begin with the situation of $H\subset G$.
Let us assume first that $G\subset \{|z|>\rho_0\}$. Let us assume that 
$\rho_G:\hat\cO_G\to G\times \hat\cO_{\{\zeta\}}$ with $\zeta$ in $H$. The map
$\rho:H\times \hat\cO_{\{\zeta\}}\to G\times \hat\cO_{\{\zeta\}}$ defined as
$\rho=\rho_G\circ \iota\circ \rho_H^{-1}$ is continuous where $\iota$ is the
inclusion. On the set
$H\times
\cO^\infty_{\{\zeta\}}$ we have $\rho(z,w)=(z,w)$. Since $H\times
\cO^\infty_{\{\zeta\}}$ is dense $\rho$ must have this form everywhere. For a set $G$
such that $f^n(G)\subset \{|z|>\rho_0\}$ we apply the function $f^n$ and repeat the
previous argument. \qed

\give Remark.  In Proposition 2.7 we have defined an atlas $\cA$ of charts which
gives us a Riemann surface lamination structure, which has the special property
that a global function,
$\varphi^+$, is used as the local holomorphic coordinate.  It is sometimes useful to
consider a related atlas
$\cA'$ which gives us an affine Riemann surface lamination.  Recall that the atlas
$\cA$ consists of charts
$\rho_G$ where
$G$ is a simply connected subset of $\C$ and 
$\rho_G:\cO_G\to G\times Y$. Define $\rho'_G=(\ell,id)\circ\rho_G$ where $\ell$ is
a branch of the logarithm function defined on $G$. The overlap functions for
$\cA'$ now have the form $(z,w)\mapsto(z+c, g(w))$ where the constant $c$ is an
integral multiple of
$2\pi$ which arises because of ambiguity in the choice of the branch of the
logarithm. The chart $\cA'$ induces an affine structure on each leaf, and with respect
to these affine structures the map $f$ has the form $z\mapsto d\cdot z+c(w)$ on the
local plaque with $w$ constant.
  Note further
that our transition functions preserve both factors of the product so that not only are
leaves well defined but local transversals to the leaves of the form $\{(z_0,y):y\in
Y_j\}$ are also well defined, independent of the chart (compare Sullivan's comment on the
TLC property [S, p. 549]). 

\bigskip
The following are consequences of Proposition 2.7.
\proclaim Corollary 2.8. Let $\cO'$ be a leaf of the lamination
$\cM^-$, then $\varphi^+|_{\cO'}$  has no critical points.

\proclaim Corollary 2.9.  Let $\cO'$ be a leaf of the lamination $\cM^-$. Then when
$\cO'$ is given the leaf topology the function $\varphi^+|_{\cO'}:
\cO'\to\cx{}-\bar\Delta$ is a covering map. 

\give Proof.  Let  $G$ be any set arising in the previous proposition. Then the
proposition implies that $\varphi^+|_{\cO'}^{-1}(G)$ consists of path components
mapped bijecitvely to $G$.  With respect to the leaf topology each of these path
components is an open set.\qed

\proclaim Lemma 2.10. Each leaf of $\cM^-$ is a conformal disk, and $\varphi^+$ is
a minimal harmonic function on each leaf.

\give Proof.  Let $\cO'$ be a leaf $\cM^-$.  Assume first that the degree of the covering
$\varphi^+|_{\cO'}$ is finite.  For $n$ a natural number let
$\cO'_n=f^{-n}(\cO')$.   Now $\varphi^+\circ
f^n=\pi^n\circ\varphi^+$ where we use the letter $\pi$ for the $d$-th power map.  The
degree of $f^n|_ {\cO'_n}$ is one and the degree of $\pi^n$ is $d^n$. Since degrees
multiply we have
$deg(\varphi^+|_{\cO'_n})\cdot d^n=deg(\varphi^+|_{\cO'})$.  Since the right-hand
side is divisible by some highest power of $d$ this equation cannot be valid for all $n$.
We conclude that the degree of
$\varphi^+|_{\cO'}$ must be infinite. It follows from covering space theory that the
image of the map
$(\varphi^+)_*:\pi_1(\cO')\to\pi_1(\cx{}-\bar\Delta)$ is a subgroup of infinite index.
But this second group is isomorphic to
$\Z$ hence the only subgroup of infinite index is the trivial group. Since by covering
space theory
$(\varphi^+)^*$ is injective we conclude that $\pi_1(\cO')$ is trivial.  Thus $\cO'$ is a
simply connected Riemann surface. 

Let $H=\{x+iy:x>0\}$. Consider the map $\alpha:H\to \C-\Delta$ given by
$\alpha(z)=e^z$. Now $\cO'$ and $H$ are both universal covering spaces of
$\C-\Delta$ so by covering space theory there is a bijection $\beta:\cO'\to H$ so
that $\varphi^+=\alpha\circ\beta$. Since $\alpha$ and $\varphi^+$ are holomorphic,
$\beta$ is holomorphic. Thus $\cO'$ is holomorphically equivalent to the right
half-plane and $\log|\varphi^+|$ is holomorphically equivalent to
$\beta(x+iy)\mapsto\log|\alpha(x+iy)|=x$. The function $x+iy\mapsto x$ is
 a minimal function on the right half-plane so $\log|\varphi^+|$ is minimal on
$U$. \qed

If we fix a simply connected domain $G$ as above and consider $\zeta_1, \zeta_2\in
G$, then the product structure on $\hat\cO_G$ gives us a homeomorphism
$\chi_{\zeta_1, \zeta_2}:\hat\cO_{\{\zeta_1\}}\to \hat\cO_{\{\zeta_2\}}$. If
$\gamma$ is a path in $\cx{}-\bar\Delta$, then we define a holonomy map
$\chi_\gamma:\hat\cO_{\{\gamma(0)\}}\to \hat\cO_{\{\gamma(1)\}}$ by
subdividing the path into small intervals, covering each interval by a disk and
composing successive local holonomy maps. The resulting holonomy map depends
only on the homotopy type of the path relative to its endpoints.

We now prove that $\hat\cO$ supports a unique positive, closed current.
 We will define a family $\cL(\cM^-)$ of positive,
closed currents which are compatible with the laminar structure of
$\cM^-$.  We consider currents
$S$ such that for each $\zeta\in\cx{}-\bar\Delta$ there is a neigborhood $G$ of
$\zeta$ and a measure
$\lambda_{S,\zeta}$ on $\hat\cO_{\{\zeta\}}$ such that
$$S\contract\hat\cO_G=\int_{t\in \hat\cO_{\{\zeta\}}}
\lambda_{S,\zeta}(t)[\Gamma_t].\eqno(2.2)$$
where $\Gamma_t$ are the leaves of the product lamination $\cM^-\cap J^-_{G}$. 
We let $\cL(\cM^-)$ denote the set of positive, closed currents $S$ on
$U^+$ with support in $\hat\cO$ such that the representation (2.2) holds with
$\lambda_{S,\zeta}$ a probability measure.  Thus for each
$S\in\cL(\cM^-)$ there is a probability measure $\lambda_\zeta$ on
$\hat\cO_{\{\zeta\}}$ for each $\zeta\in\cx{}-\bar\Delta$.  These measures are
connected via the holonomy map: if $\rho$ is a path connecting the points $\zeta$
and
$\zeta'$, then
$(\chi_\rho)_*\lambda_\zeta=\lambda_{\zeta'}$.   Let $\gamma$ be a loop based at
$\zeta$ which generates the fundamental group of $\cx{}-\bar\Delta$.  If
we have a measure $\lambda$ on
$\hat\cO_{\{\zeta\}}$ which is invariant under the holonomy automorphism
$\chi_\gamma: \hat\cO_{\{\zeta\}}\to \hat\cO_{\{\zeta\}}$ we can define  a
family of measure $\lambda_{\zeta'}$  on each
$\hat\cO_{\{\zeta'\}}$ by the formula
$\lambda_{\zeta'}=(\chi_\rho)_*\lambda_\zeta$ where $\rho$ is a path connecting
$\zeta$ and $\zeta'$.  The resulting measure is independent of the path: any other
path between these points is homotopic relative to its endpoints to a path which differs
from $\rho$ by a multiple of $\gamma$.
This family of measures produces a current in
$\cL(\cM^-)$ via formula (2.2). Thus we have proved the following Lemma.

\proclaim Lemma 2.11.  Elements of $\cL(\cM^-)$ are in one-to-one correspondence
with probability measures on $\hat\cO_{\{\zeta\}}$ which are invariant under the
holonomy map $\chi_\gamma$.

Next we will look at the push-forward of a current under a diffeomorphism.  We
recall that for a current of integration $[D]$, the push-forward is given by
$f_*[D]=[f(D)]$.  Thus the push-forward of a current of the form (2.1) is given
by
$f_*(\int\lambda(t)[\Gamma_t])=\int\lambda(t)[f\Gamma_t]$.

\proclaim Lemma 2.12.  $d^{-1}f_*\cL(\cM^-)\subset\cL(\cM^-)$.

\give Proof.  For a point $\zeta\in\cx{}-\bar\Delta$, we let
$\xi_1,\dots,\xi_d$ denote the solutions to $\xi^d=\zeta$, and we let
$D_1,\dots,D_d$ denote the preimages of a neighborhood $G$ of $\zeta$.  Let
$S\in\cL(\cM^-)$ be given, and suppose that over each $D_j$, 
$S$ has the form (2.1) .  It follows that on the 
neighborhood $\hat\cO_G$ of the fiber $\hat\cO_{\{\zeta\}}$, we have
$$f_*S\contract(\varphi^+)^{-1}G=\sum_{j=1}^d \int_{t\in
\hat\cO_{\xi_j}}\lambda_{\xi_j}(t)[f\Gamma_{j,t}],$$
where $\Gamma_{j,t}$ denotes a leaf of $\cM^-$ lying over $D_j$ and containing
the point $t$.  The transversal measure in this case is
the push-forward of $\sum_{j=1}^d\lambda_{\xi_j}$, which has total mass $d$. 
After we normalize by multiplying by $d^{-1}$, $f_*S$ will again belong to
$\cL(\cM^-)$. 
\qed

By (1.1) it follows that if  $R$ is large, then for any $S\in\cL(\cM^-)$, 
the restriction $S\contract\{|y|>R\}$ is closed on
$\{|y|>R\}$.  We will let $S_\zeta$ denote the slice measure $S|\{y=\zeta\}$.
For $|\zeta|>R$, the transversal
$\hat\cO\cap\{y=\zeta\}$ is approximately $\hat\cO_{\{\zeta\}}$, and so the local 
holonomy map $\chi:\hat\cO\cap\{y=\zeta\}\to \hat\cO_{\{\zeta\}}$ is well defined
and takes the slice measure
$S_\zeta$ to $\lambda_\zeta$.  Thus the slice measure is a probability measure.
We define the potential function
$$P_S(x,y):=\int_{x'\in\cx{}}\log|x-x'|\,S_y(x').$$
Since $S$ is a positive, closed current, it follows that on $\{|y|>R\}$, $P_S$ is
pluri-subharmonic, and ${1\over2\pi}dd^cP_S=S$.  Since each slice
measure $S_y$ is supported in the set $\{|x|\le|y|\}$, it follows that
$$\log(|x|+|y|)\ge P_S(x,y) \ge\log(|x|-|y|)\eqno(2.3)$$
holds on the set $\{(x,y):R<|y|<|x|\}$.

The following proposition was motivated by the result of Fornaess and Sibony
[FS, Theorem 7.12]: {\sl Any positive closed current on $\cx2$ which has support
in $K^+$ must be a multiple of $\mu^+$.}

\proclaim Proposition 2.13.  $\cL(\cM^-)=\{\mu^-\contract U^+\}$. 

\give Proof.  By Lemma 2.11, we know that there is an element
$S\in\cL(\cM^-)$.  If we write $S^j:=d^jf^{-j}_*S$, then by Lemma 2.12, we have
$S^j\in\cL(\cM^-)$.  By our notation, we have
$$S=d^{-j}(f^{-j})^*d^jf^{-j}_*S =d^{-j}(f^{-j})^*S^j.$$
By the remarks on the potential function, then, we have
$$S=(d^{-j}){1\over2\pi}dd^c(P_{S^j}\circ f^{-j}).$$
For any point $(x,y)\in\{|y|>R\}-K^-$, we observe that
$f^{-j}(x,y)=(x_{-j},y_{-j})$ satisfies $R\ll|y_{-j}|\ll|x_{-j}|\to\infty$ as
$j\to\infty$.  We recall that $d^{-j}\log|x_{-j}|\to G^-$ uniformly on
compact subsets of $\cx2-K^-$.  Now if we apply (2.3), we obtain that 
$$d^{-j}P_{S^j}(f^{-j}(x,y)) \to G^-(x,y)\eqno(2.4)$$
uniformly on compact subsets of  $\{|y|>R\}-K^-$.  It follows that the sequence
on the left hand side of (2.4) is uniformly bounded above, and any such
sequence of plurisubharmonic functions either has a subsequence that converges
in $L^1_{loc}(|y|>R)$ or the entire sequence converges everywhere to $-\infty$. 
Now, passing to a subsequence, we may assume that it converges in
$L^1_{loc}(|y|>R)$ to a plurisubharmonic function $W$.  Since $W=G^-$ on the
complement of $K^-$, $G^-$ is continuous and $G^-=0$ on $K^-$, and since $W$
is upper semicontinuous, we have $W=0$ on $\partial K^-$.  Further, since $\hbox{\rm
supp}(dd^cP_{S^j})\subset J^-$, it follows that $\hbox{\rm supp}(dd^cW)\subset
J^-$.  Thus $W$ is pluriharmonic on the interior of $K^-$, and so $W=0$ on
$K^-$.  Thus $W=G^-$, so this sequence converges in $L^-_{loc}(|y|>R)$ to
$G^-$.  Applying ${1\over2\pi}dd^c$ to (2.3), we conclude
that $S={1\over2\pi}dd^cG^-=\mu^-$ on $U^+$.  Thus $S=\mu^-\contract U^+$.
\qed

\proclaim Corollary 2.14.  $\bar\cO\cap U^+=\hat\cO=J^-_+$.

\give Proof. Since $\cO\subset J^-$ and $J^-$ is closed we have $\bar\cO\cap U^+\subset
J^-_+$.  For the other containment, let 
$\nu^-$ denote the current with  $\hbox{\rm supp}(\nu^-)\subset\bar\cO$, given in
Lemma 2.11.  By Proposition 2.13, we have $\mu^-\contract U^+=\nu^-$, so 
$J^-\cap U^+=\hbox{\rm supp}(\nu^-)\subset\bar\cO$, since  $\hbox{\rm
supp}(\mu^-)=J^-$.

\proclaim Corollary 2.15. The holonomy map is uniquely ergodic and the unique
invariant measure has full support.

\give Proof.  By Lemma 2.11 a finite probability measure $\lambda$ on the transversal
invariant under the holonomy map gives rise to a current on $J^-_+$.  There is only one such
current namely
$\mu^-$ and its support is all of $J^-_+$. It follows that the support of $\lambda$ is
equal to the transversal.

\proclaim Corollary 2.16. Every leaf of the lamination $\cM^-$ is dense.

\give Proof. This corresponds to the minimality of the holonomy map. The minimality
of the holonomy map is a consequence of the unique ergodicity with a measure of full
support (see [Wa, Theorem 6.17]).

The next result extablishes the uniqueness of the lamination $\cM^-$.

\proclaim Theorem 2.17.  If $f$ satisfies condition $(\dag)$ and $D$ is a connected Riemann
surface contained in
$J^-_+$, then $D$ is an open subset of a leaf of the lamination $\cM^-$.

\give Proof.  Let us suppose that $D$ is a complex disk such
that $D\cap J^-_+$ contains an open subset of $D$.  We recall that the Green
function $G^-$, restricted to the transversals $\hat\cO_{\{\zeta\}}$, is
continuous.  Thus
$\hat\cO_{\{\zeta\}}$ cannot have isolated points.  Let us consider a
neighborhood $N\subset J^-_+$ such that $\cM^-\cap N$ is the product lamination
$T\times G$, whose leaves are written $\Gamma_t$.  Let us suppose, for the sake
of contradiction, that $D$ is not contained in any leaf of $\cM^-$.  It
follows, then, that $D$ must intersect each leaf $M$ of $\cM^-$ in a
zero-dimensional set.  If $D$ intersects $\Gamma_{t_0}$ tangentially, then by
[BLS, Proposition 6.4],  the intersection of $D$ and $\Gamma_t$ is transverse
for all $t$ close to $t_0$ with
$t\ne t_0$.  

Thus we may assume that $D$ and $\Gamma_t$ intersect transversally for all 
$t\in T$.  We may also assume that in local coordinates
$D=\{|x|<1\}\times\{0\}$, so that $D$ is contained in the family of complex
disks $D_s=\{|x|<1,y=s\}$.  It follows that for $\epsilon$
small, there is a holonomy map
$\chi_s:D_s\cap N\to  D\cap N$ for each $|s|<\epsilon$.  Since $D\cap N$
contains an open subset of $D$, it follows that $D_s\cap N$ contains a nonempty
open subset of $D_s$ for each $s$.  It follows that
$\bigcup_{|s|<\epsilon}(D_s\cap N)$ contains a domain in $\cx2$, which
contradicts the fact that $J^-_+$ has no interior.  \qed

\proclaim Corollary 2.18. If $f$ satisfies condition $(\dag)$ then decomposition of
$J^-_+$ into the leaves of $\cM^-$ is the unique way of writing  $J^-_+$ as a  union of
connected Riemann surfaces.

\proclaim Corollary 2.19. If \  $W\subset J^-$ is a Riemann surface conformally
equivalent to
$\C$,  then each component of $W\cap U^+$ is a disk satisfying $(\dagger)$.

\give Remark. If $p\in\cR$ then $W=W^u(p)$ satisfies the hypotheses of this Corollary.

\give Proof. The set $U^+\cap W$ is an open subset of $W$ so it has a most countably many
components. According to Theorem 2.17 each component is contained in a leaf $L$ of the
lamination $\cM^-$. Let $L_1,L_2,\ldots$ be the leaves that meet $W$. Let
$X=W\cup\{L_j\}$. We claim that $X$ has the structure of an immersed Riemann surface. 
Let $\phi:\C\to W$ and let $\phi_j:H\to L_j$ where $H$ is the right half plane. The set of
charts
$\cA=\{\phi,\phi_1,\ldots\}$ forms an atlas for $X$. 
 To show that $\cA$ is an atlas we verify that the change of coordinate functions are
holomorphic. If $\phi(\C)\cap\phi_j(H)=W\cap L_j$ is a Riemann surface so $\phi^{-1}(W\cap
L_j)$ is an open subset of $\C$ and $\phi_j^{-1}(W\cap L_j)$ is an open subset of $H$ and
$\phi^{-1}\circ\phi_j$ is holomorphic.

The atlas $\cA$ determines a topology on $X$ where we define a set $Y\subset X$ to be open if
$\phi^{-1}(Y)$ is open for each chart $\phi\in \cA$. With respect to this topology the
inclusion from $X$ into $\C^2$ is continuous. The existence of an atlas for a toplogical
space does not imply that it is Hausdorff. From the fact that the inclusion is
continuous we can deduce that
$X$ is Hausdorff. If $p$ and $q$ are distinct points in $X$ then they are contained in disjoint
neighborhoods in $\C^2$. The pullbacks of these neighborhoods are disjoint neighborhoods in
$X$. We also see that $X$ is path connected. Since $X$ is holomorphically embedded in
$\C^2$ so it is not compact.

We will finish by showing that $X=W$. This will imply that the intersection of $W$ with $U^+$
is a union of leaves of the lamination $\cM^-$. The corollary then follows from
Theorem 2.1 (5).

Assume that $W$ is a proper subset of $X$. Replace $X$ by its universal cover
$\tilde X$ and replace $W$ by some lift $W_0$ of $W$. The inclusion of $W_0$ in
$\tilde X$ is still proper. Now $\tilde X$ is simply connected but not compact
thus it is conformally equivalent to a disk or plane.  Now $\tilde X$ contains
a copy of $\C$ so it is not the disk. We conclude that $\tilde X$ is the plane.
A proper simply connected subsurface of the plane is equivalent to the disk so
$W_0$ must be equal to $\tilde X$ contradicting our assumption. \qed

\section 3. External Rays and Semi-conjugacy to the Solenoid 

The previous section was devoted to the study of some of the properties of mappings which
have a disk satisfying (\dag), and the following sections will be devoted to establishing
various conditions which imply this property.  Before we proceed with
this, however, we describe two ways in which the special features of these maps can be
exploited.  First we define a family $\cE$ of external rays in $J^-_+$.  It is
expected that, as in the 1-dimensional case, a connection between the dynamics of
the restriction $f|_{J^-_+}$ and $f|_J$ will be obtained by a more careful study
of the map of external rays.  Our second observation is that 
$f|_{J^-_+}$ is semi-conjugate to the shift on the complex solenoid.  The paper [BS7] is
devoted to further exploration of the relationship between these mappings and the solenoid in
the hyperbolic case.

We begin by recalling some properties developed in \S2.  If $f$ satisfies
(\dag), then by Theorem 2.1 there is a lamination $\cM^-$ of $J^-_+$ and a
holomorphic extension $\varphi^+:J^-_+\to\cx{}-\bar\Delta$, which serves as a
canonical projection.  For any $S\subset\cx{}-\bar\Delta$ we set
$J^-_S:=(\varphi^+)^{-1}(S)$.  The sets $J^-_{\{\xi\}}$,
$\xi\in\cx{}-\bar\Delta$ form a canonical family of transversals to $J^-_+$.  If
$G\subset\cx{}-\bar\Delta$ is a simply connected domain, then for any $\xi\in
G$, Proposition 2.4 shows that the restriction $\cM^-\cap J^-_G$ is equivalent to a
product lamination whose leaves are $\{\Gamma_t:t\in
J^-_{\{\xi\}}\}$, and 
$\varphi^+:\Gamma_t\to G$ is a canonical biholomorphism.  This local  product
structure also extends to $\mu^-$ in the sense that there is a measure
$\lambda_\xi$ on $J^-_{\{\xi\}}$ such that 
$$\mu^-\contract J^-_G=\int_{t\in
J^-_{\{\xi\}}}\lambda_\xi(t)[\Gamma_t].\eqno(3.1)$$

For a leaf $M$ of $\cM^-$, the 1-form $d^cG^+$ restricts naturally to $M$. 
The set of integral curves $\gamma$ of $d^cG^+|_M$ will be called {\it external
rays}, and the set of all external rays will be denoted by $\cE$.  External rays
serve as a substitute for gradient lines of $G^+$, and each $\gamma$  may
naturally be parametrized by $G^+$.  These are called rays because each
$\gamma\in\cE$ is contained in some leaf $M$ of $\cM^-$, and $\gamma$ is the
lift of some radial ray $R_\theta=\{re^{i\theta}\in\cx{}:r>1\}$ under the
mapping $\varphi^+|_M$.  We may define the map $e_r:\cE\to\{G^+=r\}\cap
J^-$ by letting $e_r(\gamma)$ be intersection point  $\gamma\cap\{G^+=r\}$.  The mapping 
$\chi_{r,s}:J^-\cap\{G^+=r\}\to J^-\cap\{G^+=s\}$, defined by following the
external rays, is a homeomorphism between these two sets, and
$e_s=\chi_{r,s}\circ e_r$.  It is natural to topologize
$\cE$ so that $e_r$ is a homeomorphism. 

We will show how external rays give a description of the
harmonic measure $\mu$ when (\dag) holds. Let $\mu_r$ denote the measure on
$J^-\cap\{G^+=r\}$ defined by 
$$\mu_r:={1\over (2\pi)^{2}}(-d^cG^+\contract\partial\{G^+>r\})\wedge dd^cG^-
={1\over (2\pi)^{2}}dd^c(\max(G^+,r))\wedge dd^cG^-.\eqno(3.2)$$  
We note that $(2\pi)^{-1}(-d^cG^+|_M\contract\{G^+>r\})$ may be interpreted as the
harmonic measure of
$M\cap\{G^+>r\}$. (This is also the pullback under $\varphi^+|_M$ of the
planar harmonic measure $(2\pi)^{-1}d^c\log|\zeta|$ of $\{\zeta\in\cx{}:|\zeta|>r\}$,
which is normalized arclength measure on the circle
orthogonal to the rays.)  For a simply connected domain
$G\subset\cx{}-\bar\Delta$, it follows from (3.1) that on $J^-_G$ the middle
expression in (3.2) may be interpreted as
$$\int_{t\in J^-_{\{\xi\}}}\lambda_\xi(t)\omega_t \eqno(3.3)$$
if we set $\omega_t:=(2\pi)^{-1}(-d^cG^+|_{\Gamma_t}\contract\{G^+>r\})$.  Flowing
along the gradient lines of $G^+|_M$ preserves harmonic measure in the sense that
${\chi_{r,s}}_*\omega_r=\omega_s$, wherever $\chi_{r,s}$ is correctly defined, so
$\chi_{r,s}$ also preserves
$\mu_r$ in the sense that 
$(\chi_{r,s})_*\mu_r=\mu_s$.  Thus we may define a
measure $\nu$ on $\cE$ by the condition that $({e_r})_*\nu=\mu_r$, and this definition
is independent of $r$.

For a leaf $M$ of $\cM^-$, let $\alpha_M:H\to M$ denote the uniformization given
in Lemma 2.2.  All external rays $\gamma\subset M$ are of the form 
$\gamma(y)=\{\alpha_M(x+iy):0<x<\infty\}$.  Since $\alpha_M$ is a bounded
analytic function on the set $\{0<x<R<\infty\}$ the limit
$\alpha_M^*(y):=\lim_{x\to0^+}\alpha_M(x+iy)$ exists for almost every $y$. 
Thus $\alpha_M^*(y)=\lim_{r\to0^+}e_r(\gamma(y))$ exists for almost every
$\gamma(y)\subset M$.  Since harmonic measure $\omega_t$ corresponds to $dy$ under
$\alpha_M$, it follows that within each leaf $M$, $e$ is defined almost
everywhere with respect to $\omega_t$.  We define the endpoint map
$$e(\gamma)=\lim_{r\to0}e_r(\gamma)$$
for all $\gamma\in\cE$ for which this limit exists, and we observe that by
(3.3),
$e(\gamma)$ is defined for $\nu$ almost every $\gamma$.  

Thus $e:\cE\to\cx2$ is a Borel measurable map, and we may push the measure
$\nu$ forward to a measure $e_*\nu$ on $\cx2$, which coincides with
the weak-$*$ limit $\lim_{r\to0}(e_r)_*\nu$.  From the right hand expression in
(3.2) and the fact that $G^+$ is continuous, we see that $\lim_{r\to0^+}\mu_r=\mu$. It
follows by the uniqueness of weak-$*$ limits of measures that we have the following:

\proclaim Theorem 3.1.  If $f$ satisfies the hypotheses of Theorem 2.1, then the
endpoint mapping $e$ is defined $\nu$ almost everywhere on $\cE$, and
$e_*\nu=\mu$.

\bigskip

We now begin the construction of the semi-conjugacy to the standard model.
Let $\sigma:\cx*\to\cx*$ denote the map $\sigma(z)=z^d$.  Let $\Sigma$ denote
the projective limit of the map $\sigma$, and call $\Sigma$ the {\it $d$-fold
complex solenoid}.  We may represent the solenoid as a space of bi-infinite
sequences
$$\Sigma=\{(\dots
z_{-2}z_{-1}z_0z_1z_2\dots):\hbox{\rm where }
z_j\in\cx*\hbox{\rm and }z_{j+1}=\sigma(z_j) \hbox{ \rm  for all } j\in\ZZ\}.$$
Observe that $\sigma$ induces a homeomorphism (which we again denote as
$\sigma:\Sigma\to\Sigma$), which is given by left translation: $\sigma(z)=w$,
where $z=(z_j)$ and $w=(w_j)$ has entries $w_j=z_{j+1}$.

Define $\pi:\Sigma\to\cx*$ by $\pi(z)=\pi(\dots z_{-1}z_0z_1\dots)=z_0$.  Thus
$\pi\sigma=\pi^d$.  We give $\Sigma$ the product topology, so the fiber
$\pi^{-1}(\zeta)$ is a Cantor set (i.e.\ totally disconnected and perfect).

The standard (real) solenoid is given by 
$$\Sigma_0=\{s\in\Sigma:|\pi(s)|=1\}.$$
This may also be identified with a set of sequences of points of the circle
$$\Sigma_0=\{\theta=(\dots\theta_{-1}\theta_0\theta_1\dots):
\theta_j\in\R/2\pi\ZZ, \  d\cdot\theta_n=\theta_{n+1}\}.$$
We let $\eta:\Sigma\to\Sigma_0$ be defined by $\eta(s)=\tilde s$, where $\tilde
s_j=s_j/|s_j|$.  Thus $\eta$ commutes with $\sigma$.  We define
$$\Sigma_+=\{s\in\Sigma:|\pi(s)|>1\}$$
to be the portion of the complex solenoid lying above the complement of the
closed unit disk.  We observe that the fibers of the mapping
$\eta:\Sigma_+\to\Sigma_0$ are rays in the complex solenoid, so the complex
solenoid has a natural structure of a family of external rays, parametrized by
the real solenoid $\Sigma_0$.

Given a space $X$, a mapping $\Phi:X\to\Sigma$ is given by a family of mappings
$\phi_j:X\to\cx*$,
$j\in\ZZ$ such that $\phi_{j+1}(x)=\phi_j^d(x)$ for all $x\in X$ and
$j\in\ZZ$.  If $\Phi$ is to give a semi-conjugacy between a bijection $f:X\to X$
and $\sigma$, then we must have $\phi_j\circ f=\phi_{j+1}=\phi_j^d$.  In
this case $\phi_0:X\to\cx*$ determines all of the coordinate maps $\phi_j$ via the relation
$\phi_j(x)=\phi_0(f^jx)$ for all $j\in\ZZ$.  The consistency condition for a
given map $\phi_0:X\to\cx*$ to induce an equivariant mapping $\Phi$ in this fashion is that
$\phi_0\circ f=\phi_0^d$.  According to Theorem 2.1 (1) if there is a disk satisfying
condition $(\dagger)$ then 
$\varphi^+$ has an extension to
$J^-_+$ which satisfies 
$\varphi^+\circ f=(\varphi^+)^d$. Thus $\varphi^+$ serves as a $0$-th coordinate map for an
equivariant  mapping
$\Phi:J^-_+\to\Sigma_+$ given by
$$\Phi(p):=(\phi_j(p))=(\varphi^+(f^jp)).\eqno(3.4)$$

\proclaim Theorem 3.2. If $f$ satisfies the hypotheses of Theorem 2.1, then
there is a continuous mapping $\Phi:J^-_+\to\Sigma_+$ which is holomorphic and
injective on the leaves of $\cM^-$ and such that
$$\sigma\circ\Phi=\Phi\circ f,\eqno(3.5)$$
and
$$\log|\pi\circ\Phi|=G^+.\eqno(3.6)$$

\give Proof.  By Theorem 2.1, $\varphi^+$ extends to $J^-_+$ and is holomorphic on
the leaves of $\cM^-$.  Thus, if we define
$\Phi$ as in (3.4), then (3.5) and (3.6) are easily seen to be satisfied.  For
$M\in\cM^-$ it is evident that $\Phi|_M$ is holomorphic.  We must show that $\Phi|_M$
is injective.  We let
$\zeta_0=\varphi^+(p)$ and let $j\tau$ denote the path in $\cx{}$ starting at
$\zeta_0$ which traces the circle $\{|\zeta|=|\zeta_0|\}$ $j$ times in the
counter-clockwise direction.  We let
$p_j$ be the point obtained by lifting $j\tau$ via the map $\varphi^+|_M$
starting at $p_0$, so that $(\varphi^+)^{-1}(\zeta_0)=\{p_j:j\in\ZZ\}$.  We let
$\gamma_{j_1,j_2}$ denote the curve above $\tau$ which goes from $p_{j_1}$ to
$p_{j_2}$.  Since $(\varphi^+)^d=\varphi^+\circ f$, we have
$$\left((\varphi^+\circ f^{-n})|_M\right)^{d^n} =\varphi^+|M.$$
It follows that the curve $(\varphi^+\circ f^{-n})(\gamma_{j_1,j_2})$ starts at
the point $\varphi^+\circ f^{-n}(p_{j_1})$ and moves inside the circle
$\{|\zeta|=|\varphi^+\circ f^{-n}(p_0)|\}$ through an angle of
$2\pi(j_2-j_1)d^{-n}\in\ZZ$.  It follows that if
$\phi_{-n}(p_{j_1})=\phi_{-n}(p_{j_2})$ for $n=1,2,\dots$, then $j_2-j_1$ is
divisible by $d^n$ for all $n$.  Thus $j_1=j_2$, and so $\Phi$ is one-to-one on
$M$.  \qed

Under the mapping $\Phi$, the external rays $\cE$ are mapped
to rays of the solenoid, and the lamination $\cM^-$ is taken to the lamination
of $\Sigma_+$.

\section 4.  Unstable Connectedness

Let $p$ be any point in $J$ so that $W^u(p)$ exists and is conformally equivalent to $\C$. 
We will say that $f$ is unstably connected with respect to $p$ if $W^u(p)-K^+$ has at least
one simply connected component.

 Let $\nu$ be an ergodic hyperbolic measure with index one. 
We will say that $f$ is {\it unstably connected} with respect to $\nu$ if for $\nu$ almost
every point $p$, $f$ is unstably connected with respect to $p$.   Since $\nu$ is
assumed to be ergodic and our condition on $p$ is invariant under $f$, it is equivalent to
assume $f$ is unstably connected with respect to $p$ for $p$ in a set of positive $\nu$
measure. 

\proclaim Theorem 4.1. Let $\nu$ be a hyperbolic index one measure then $f$ is
unstably connected with respect to $\nu$ if and only if there is a disk satisfying
(\dag).

\proclaim Corollary 4.2. If $f$ is unstably connected with respect to some
hyperbolic measure
$\nu$ then it is unstably connected with respect to any hyperbolic measure $\nu$.

We say that $f$ is {\it unstably connected} if $f$ satisfies the hypothesis of the corollary
for any hyperbolic measure.

\give Remark.  If $p$ is a periodic saddle point and $\nu$ is the normalized counting measure on the orbit
then the condition of being unstably connected with respect to $p$ described in the
introduction, i.e. that $W^u(p)\cap U^+$ has a simply
connected component, is equivalent to being unstably connected with respect to $\nu$.

The ``if'' part of Theorem 4.1 is contained already in Corollary 2.19. 
The rest of this section will be devoted to the proof of the ``only if'' satement. 
The proof of Theorem 4.1 can be broken into two parts: function
theoretic, involving the growth of $G^+$ on a single unstable manifold, and 
dynamical, where we use the properties of the measure $\nu$.   We
will start with some properties of hyperbolic measures.  The reader who is
unfamiliar with hyperbolic measures may find it convenient upon first reading to
consider only the simplest hyperbolic measures: the averages of
point masses over periodic saddle orbits.  This proves Theorems 4.1 in the case of
periodic saddle points and makes the proofs of the following three Lemmas 
rather simple.

For the rest of this Section, $\nu$ will denote an index one
hyperbolic measure, and $f$ will be assumed to be unstably connected with
respect to $\nu$.  

Applying the Pesin Stable Manifold Theorem to an index one hyperbolic measure
$\nu$, it follows that for $\nu$ almost every point
$p$ there is an unstable manifold $W^u(p)$ which is conformally equivalent to
$\cx{}$ (cf.\ [BLS], [W]).  Let us fix a conformal equivalence
$\phi:\cx{}\to W^u(p)$ with $\phi(0)=p$.  We will use $\phi$ to translate
concepts from $W^u(p)$ back to $\cx{}$.  For $z\in\cx{}$ we will sometimes write
$G^+(z)$ to denote the function $G^+(\phi(z))$.  We write $U^+$ for $\phi^{-1}(U^+)$
and
$K^+$ for
$\phi^{-1}(K^+)$.  Thus $G^+$ is a continuous subharmonic function on $\cx{}$
which is equal to zero on $K^+$ and positive and harmonic on $U^+$.

Given one uniformization $\phi$, all other uniformizations $\phi$ are of
the form $\phi(\alpha z)$ for some constant $\alpha\in\cx{}-\{0\}$.  We define
the function
$$M(p,r)=\max_{|z|=r}G^+(\phi(z)).$$

By the maximum principle $M$ is an increasing function of $r$.
We may choose the scale $|\alpha|$ in the uniformizing function $\phi(\alpha z)$
so that $M(p,1)=1$, and we let $\phi_p$ denote the uniformization which is
normalized this way.  This determines
$M$ uniquely in terms of
$G^+$ and the conformal structure of $W^u(p)$. Since $G^+\circ\phi$ is subharmonic on
$\cx{}$,
$M(p,r)$ is continuous and increasing in $r$.  We define a Hermitian metric
$\Vert\cdot\Vert_G$ on the subspace $E^u_p$ by the condition
$$\Vert D\phi_p(0){\bf 1}\Vert_G =1$$
where ${\bf 1}$ denotes an element of $T_0\cx{}$ which has unit length with
respect to the standard metric.  We will also consider the growth rate of
$Df^n$ with respect to $F$:
$$\Vert Df^n_p|_{E^{u}_p}\Vert_G=\Vert Df^n_p(v)\Vert_G/\Vert
v\Vert_G,$$
where $v$ is any nonzero element of $E^u_p$.  The advantage of
$\Vert\cdot\Vert_G$ is that it transforms naturally under $f$ in
connection with $G^+$ and $M(p,r)$.

\proclaim Lemma 4.3.  For all $n\in\ZZ$ we have $M(f^np,\Vert
Df^n_p\Vert_G)=d^n$.

\give Proof.  The map $\phi^{-1}_{f^np}\circ f^n\circ\phi_p$ is a holomorphic
and bijective map that sends $\cx{}$ to $\cx{}$ and takes 0 to 0.  Any such map
is linear.  Thus there is $\xi\in\cx{}-\{0\}$ so that
$$\phi^{-1}_{f^np}\circ f^n\circ\phi_p(z)=\xi z.\eqno(4.1)$$
Thus $\phi_{f^np}(\xi z)=f^n\circ\phi_p(z)$.  Rewriting this equation and
applying $G^+$ gives:
$$G^+(\phi_{f^np}(\xi z))=G^+(f^n\circ\phi_p(z)).$$
Now $G^+$ multiplies by $d$  when $f$ is applied so 
$$G^+(\phi_{f^np}(\xi z))=d^nG^+(z).$$
Now we evaluate:
$$\eqalign{M(f^np,|\xi|)&=\max_{|z|=|\xi|}G^+(\phi_{f^np}(z))
=\max_{|z|=1}G^+(\phi_{f^np}(\xi z))\cr 
&=d^n\max_{|z|=1}G^+(\phi_p(z))
=d^n.\cr}$$ To evaluate $|\xi|$, we differentiate equation (4.1) to get
$$D\phi^{-1}_{f^np}\circ Df^n\circ D\phi_p(z)=\xi$$
so
$$\Vert D\phi^{-1}_{f^np}\Vert\cdot\Vert Df^n\Vert_G\cdot\Vert
D\phi_p(z)\Vert=\Vert\xi\Vert.$$
Since $\phi_p$ and $\phi_{f^np}$ were normalized so that $D\phi_p$ and
$D\phi_{f^np}$ have norm one, we have $\Vert Df^n\Vert_G=|\xi|$. \qed

The next Lemma shows that we may also compute the Lyapunov exponent starting with the
metric $\Vert\cdot\Vert_G$. 

\proclaim Lemma 4.4.  Let $\nu$ be an index one hyperbolic measure.  For $\nu$
almost every $p$, we have the existence of the limit
$$\lim_{n\to\pm\infty}{1\over n}\log\Vert D^nf_p\Vert_G=\lambda^u(\nu).$$

\give Proof.  Let $v\in E^u_p$ be nonzero, and let
$$r(p)={\Vert v\Vert_H\over\Vert v\Vert_G},$$
so that we have
$${\Vert Df_p\Vert_H\over\Vert Df_p\Vert_G} = {\Vert Df_p(v)\Vert_H\over\Vert
Df_p(v)\Vert_G} \cdot {\Vert v\Vert_G\over\Vert v\Vert_H} = r(f(p))/r(p).$$
Let $\alpha(p)=\log\Vert Df_p\Vert_H$,   $\beta(p)=\log\Vert Df_p\Vert_G$,
and   $\rho(p)=\log r(p)$.  Taking logarithms gives
the cocycle equation:
$$\alpha(p)-\beta(p)=\rho(f(p))-\rho(p).$$  If $\rho$ were in $L^1$, then our
Lemma would be a consequence of the Ergodic Theorem.  

Although we have no information on $\rho$, we can show that 
$\alpha(p)-\beta(p)$ is bounded below.  Since $M(fp,||Df_p||_G)=d$ and
$M(fp,1)=1$ and $M$ is strictly monotone in $r$, we conclude that
$\Vert Df_p\Vert_G>1$.  Thus
$\alpha(p)>0$.  Further, $\Vert Df_p\Vert_H$ is bounded above by the supremum
of the Euclidean norm of the Jacobian matrix $Df_p$ over the compact set $J$. 
Thus $\beta(p)\le C$ and $\alpha(p)-\beta(p)\ge-C$.  According to [LS,
Proposition 2.2]  we have $\lim_{n\to\infty}{1\over
n}(\alpha(f^np)-\beta(f^np))=0$ so
$\lim_{n\to\infty}{1\over n}\alpha(f^np) = \lim_{n\to\infty}{1\over
n}\beta(f^np)$ holds for $\nu$ almost every point $p$, and this limit must be
equal to $\lambda^u(\nu)$.  \qed

We fix a component $\cO$ of $U^+\cap W^u(p)$.  An {\it end} of $\cO$, written
$\cE_r(\cO)$ or just $\cE_r$, is a connected component of $\cO\cap\{|z|>r\}$.  It
is evident that $\cE_0(\cO)=\cO$.  We say that an end $\cE_r$ has {no loops}
if $\cE_r$ contains no closed curves which encircle the origin. 

\proclaim Lemma 4.5.   $W^u(p)\cap U^+$ has no loops.

\give Proof.  If $U^+$ contains a loop around zero then the component of $K^+$
containing $0$ is compact. But we have assumed that there are no compact
components. \qed

For an end $\cE_r\subset U^+_p$, we define
$$M(p,\cE_r,s)=\max_{\zeta\in\bar\cE_r\cap\cO\cap\{|\zeta|=s\}}G^+(\zeta),$$  
or we write $M(\cE_r,s)=M(p,\cE_r,s)$ if the point $p$ is understood.  
Let us fix $r<s$.  For a subset $X\subset\partial(\cO\cap\{r<|z|<s\})$, we
define a function $\omega(X,z)$ called the harmonic measure of the set $X$ to be 
to be the greatest function
$0\le\omega(X,z)<1$ that is harmonic on $\cO\cap\{r<|z|<s\}$ and which satisfies 
$$\lim_{\zeta\to\zeta_0}\omega(X,\zeta)=0\hbox{ 
\ \ for }\zeta_0\in\bd(\cO\cap\{r<|z|<s\})-X.$$

We will consider the case where $X$ has the form $A_t=\cO\cap\{|z|=t\}$ for
$t=r,s$.  Let $t\Theta(t)$
denote the length of the set $\cO\cap\{|z|=t\}$.  When $\cO$ contains no
loops, classical estimates on the harmonic measures of the sets $A_r$ and $A_s$ for
$r<|z|<s$ (see Fuchs [F]) are given by
$$\omega(A_r,z)\le 4\cdot\hbox{\rm exp}\left(-\pi\int_r^{|z|}{dt\over
t\Theta(t)}\right)$$
and
$$\omega(A_s,z)\le4\cdot\hbox{\rm exp}\left(-\pi\int_{|z|}^s{dt\over
t\Theta(t)}\right).$$
Since $\Theta(t)\le2\pi$ we have
$$\pi\int_a^b{dt\over t\Theta(t)}\ge \pi\int_a^b{dt\over 2\pi t}
={1\over2}\int_a^b{dt\over t}={\log(b/a)\over 2}.$$
Thus 
$$\omega(A_r,z)\le4\sqrt{r/|z|}\quad{\rm\ and\
}\quad\omega(A_s,z)\le4\sqrt{|z|/s}.$$

We use these estimates as follows.  For any end $\cE_r$ of $\cO$, the maximum
principle gives $$G^+(z)\le M(\cE_r,r)\omega(A_r,z) + M(\cE_r,s)\omega(A_s,z).\eqno(4.2)$$
Thus
$$M(\cE_r,|z|)\le 4M(\cE_r,r)\sqrt{r/|z|}+4M(\cE_r,s)\sqrt{|z|/s}.\eqno(4.3)$$

\proclaim Proposition 4.6. Fix $r\ge0$ and an end $\cE_r$.  Then exactly one of
the following two alternatives holds.  Either
 $$M(\cE_r,s)\le 4M(\cE_r,r)\sqrt{r/s}\hbox{ for all }s>r,\eqno(4.4)$$
 or for some constant $c>0$
$$M(\cE_r,s)\ge c\sqrt s\qquad{\rm\ for\ all\ }
s>r.\eqno(4.5)$$
When $r=0$, the second alternative must hold.

\give Proof.  Let us fix an end $\cE_r$ of $\cO$.  If alternative (4.4) does not
hold, then there exists $t_0>r$ such that 
$$M(\cE_r,t_0)> 4M(\cE_r,r)\sqrt{r/t_0}.\eqno(4.6)$$
Set $\alpha=M(\cE_r,t_0)-4M(\cE_r,r)\sqrt{r/t_0}$.  Then by (4.6), $\alpha>0$
and using (4.3) with $|z|=t_0$ we have
$$\alpha\le4M(\cE_r,s)\sqrt{t_0\over s},$$
so
$$c\sqrt s\le M(\cE_r,s)$$
 holds with $c=\alpha/(4\sqrt{t_0})$.

Now in case $r=0$ we have $\cE_r=\cO$.  Case (4.4) cannot hold, for it implies
that $G^+=0$ on $\cO$ by the maximum principle.
\qed

Let $\cE_r$ be an end. If  case (4.4) in Proposition 4.6 applies, we will say that
$\cE_r$ is a decay end, and if (4.5) applies we say that, $\cE_r$ is a growth end.
By Proposition 4.6 each end is either a growth end or a decay end.
For $f>0$ we let $\cC(r)$ denote the set of connected components of
$U^+\cap\{|z|>r\}$.  Let $c(r)$ be the cardinality of the set $\cC(r)$, and let
$g(r)$ denote the number of growth components in $\cC(r)$.  Thus $c(r)\ge g(r)$. 
Since $\cC(0)$ corresponds to the set of connected components of $U^+$, it
follows from Proposition 4.6 and the maximum principle that each 
component of $U^+$ is a growth component.  Thus $c(0)=g(0)$.

If $r<s$, then there is the containment map
$$\Upsilon:\cC(s)\to\cC(r)$$
where for any component $\cO\in\cC(s)$, $\Upsilon(\cO)$ denotes the element of
$\cC(r)$ containing $\cO$.

\proclaim Proposition 4.7.  
The functions $c(r)$ and $g(r)$ are nondecreasing in
$r$, and if $g(r)$ is finite for some value of $r$, then
$\lim_{r\to0}g(r)=g(0)=c(0)$.

\give Proof.  The mapping $\Upsilon$ is  surjective since there are no bounded
components (this is a consequence of the maximum principle), so
$c(r)\ge c(s)$ for $r\le s$.  Further, by Proposition 4.6, an end $\cE_r$ is a
growth end if and only if one of the components of $\cE_r\cap\{|z|>s\}$ is a
growth end.  Thus $\Upsilon$ is a surjective mapping from growth ends to growth
ends, and so $g(r)$ is nondecreasing.  

Finally, we show that $\lim_{r\to0}g(r)=g(0)$ if $g(r)$ is finite.  Let us
suppose, to the contrary, that $g(0)<g(r)$ for $0<r<_0$.  This means that
there are two growth ends $\cE_{r_0}^1$ and $\cE_{r_0}^2$ in $U^+\cap\{|z|>r\}$
which are contained in different connected components of $U^+\cap\{|z|>r\}$ for
$0<r<r_0$, but which are contained in the same connected component of
$U^+-\{0\}$.  Now if $\gamma$ is a path in $U^+-\{0\}$ which connects these two
components, then $\gamma$ avoids some disk $\{|z|<\epsilon\}$.  Thus
$\cE_{r_0}^1$ and $\cE_{r_0}^2$ are contained in the same connected component of
$U^+\cap\{|z|>s\}$, which is a contradiction. \qed

\proclaim Lemma 4.8.  If $\cO_1,\dots,\cO_N$ are disjoint, open sets, then
there exists $j$ such that
$$\pi\int_t^s{dr\over r\Theta_j(r)}\ge{N\over 2}\log\left({s\over t}\right).$$

\give Proof.  By the Cauchy inequality we have
$$\sum\Theta_j\sum\Theta_j^{-1}\ge
\left(\sum\Theta_j^{1\over2}\Theta_j^{-{1\over2}}\right)^2=N^2.$$
Since the sets $\cO_j$ are disjoint, $\sum\Theta_j\le2\pi$, so
$$\sum_{j=1}^N{1\over\Theta_j(r)}\ge {N^2\over2\pi}.$$
This gives
$$\sum_{j=1}^N\int_t^s{dr\over r\Theta_j(r)}\ge {N^2\over2\pi}\log\left({s\over
t}\right)$$
from which the Lemma follows. \qed

\proclaim Proposition 4.9.  If for some $r>0$ there are $N$ distinct growth
ends in $U^+\cap\{|z|>r\}$, then 
$$M(U^+,s)\ge Cs^{N/2}$$
for some $C>0$ and all $s\ge r$.

\give Proof.  The harmonic measure estimates above yield the
estimate
$$M(\cO,t)\le 4M(\cO,r)\sqrt{r/t} +4M(\cO,s)\cdot\hbox{\rm
exp}\left(-\pi\int_t^s{dr\over r\Theta(r)}\right).$$
So arguing as in Proposition 4.6 we have 
$$M(\cO,s)\ge C\cdot\hbox{\rm exp}\left(-\pi\int_t^s{dr\over
r\Theta(r)}\right).$$
If we choose $\cO=\cO_j$ satisfying the conclusion of Lemma 4.8, then the
Proposition follows. \qed

\proclaim Proposition 4.10.  For $\nu$ almost every point $p$ and every $r>0$ we
have
$$g(r)\le{2\log d\over\lambda(\nu)}.$$

\give Proof.  According to Lemma 4.4, $M(f^np,\Vert Df^n_p\Vert_G)=d^n$. 
Thus, after changing variables, we have $M(p,\Vert
Df^n_{f^{-n}p}\Vert_G)=d^n$.

Let $N$ be a (finite) integer no greater than $g(r)$.  Then
$$d^n=\gamma(p,\Vert Df^n_{f^{-n}p}\Vert_G)\ge c\Vert
Df^n_{f^{-n}p}\Vert_G^{N/2}$$
and taking logarithms gives
$$n\log d\ge \log C +{N\over2}\log \Vert Df^n_{f^{-n}p}\Vert_G.$$
We divide by $n$, take limits, and apply the chain rule $\Vert
Df^n_{f^{-n}p}\Vert=\Vert Df^{-n}_p\Vert^{-1}$, so
$$ \log d\ge {N\over2}\lim_{n\to\infty}{1\over n}\log\Vert
Df^n_{f^{-n}p}\Vert_G= -{N\over2}\lim_{n\to\infty}{1\over n}\log\Vert
Df^{-n}p\Vert_G.$$
By Lemma 4.9 we have
$$\log d\ge{N\over2}\lambda(\nu).$$
Since this holds for all finite $N\le g(r)$ we have:
$$\log d\ge{g(r)\over2}\lambda(\nu).$$\qed

We say that a component $\cO$ has a unique growth end if for each $r\ge0$ there
is only one growth end $\cE_r(\cO)$.

\proclaim Theorem 4.11.  For $\nu$ almost every point $p$ every component of
$W^u_p\cap U^+$ has a unique growth end.  Furthermore the number of components
of $W^u_p\cap U^+$ is equal to a constant $N$ for $\nu$-almost every $p$, and
$N\le{2\log d\over \lambda}$.

\give Proof.  We denote by $g(p,r)$ the number of growth ends in
$\phi^{-1}_p(U^+)\cap\{|z|>r\}$.  We observe that by the property of
$\Vert\cdot\Vert_G$, we have $g(f(p),\Vert Df_p\Vert_G\cdot r)=g(p,r)$.  We
begin by showing that for almost every $p$ the function $g(p,r)$ is independent
of $r$ for $r>0$.  Since $g(f(p),1)=g(p,\Vert Df^{-1}_p\Vert_G)\le g(p,1)$, we
have $g(p,1)-g(f(p),1)\ge0$.  On the other hand
$$\int g(p,1)\nu(p) = \int g(f(p),1)\nu(p)$$
because of the invariance of $\nu$.  Both integrals are finite because $g(p,1)$
is uniformly bounded a.e.\ by Proposition 4.10.  This gives
$$\int\left(g(p,1)-g(f(p),1)\right)\nu(p)=0,$$
and so we conclude that $g(p,1)-g(f(p),1)=0$ outside a set measure
zero.  Thus
$g(f^np,1)$ is independent of $n\in\ZZ$, except for $p$ in a set of measure
zero.  Let $r_n=\Vert Df^n_{f^{-n}p}\Vert_G$, then $g(f^np,1)=g(p,r_n)$.  Since
$g(p,r)$ is a monotone function of $r$, and since $g(p,r_n)=g(p,r_m)$,
it follows that $g(p,r)$ is constant for $r_n\le r\le r_m$.  Finally,
since $r_n\to0$ as $n\to-\infty$ and $r_n\to+\infty$ as $n\to+\infty$, it
follows that $g(p,r)$ is constant on the interval $0<r<\infty$.  

Now we apply 
Proposition 4.7 to conclude that  $g(p,0)=\lim_{r\to0}g(p,r)$.  As was observed
before Proposition 4.7, $g(0)=c(0)$ is the number of connected components of
$\phi^{-1}_p(U^+)$, so the Theorem follows from Proposition 4.10. \qed

\give Remark.  In the special case where $\nu$ is normalized counting measure on
the orbit of a periodic saddle, then it is possible to define a combinatorial
rotation number in terms the action of $f$ on the growth ends.  The
Pommerenke-Levin-Yoccoz inequality gives the estimate above, as well as
information on the possible combinatorial rotation number.  Indeed, the proof of
this inequality, as given in [H], goes through essentially without change in the
case of polynomial diffeomorphisms.

\proclaim Theorem 4.12.  Let $\cO$ be a simply connected component of $W^u_p\cap
U^+$ with a unique growth end.  Then $\cO$ satisfies (\dag).

\give Proof.  Since $p$ is a Pesin regular point, it follows (see [BLS]) that
$W^u(p)\subset J^-$ and that 
$W^u(p)\cap W^u(f^np)=\emptyset$ or $W^u(p)=W^u(f^np)$ for $n\in\ZZ$.  We must
show that $G^+|_{\cO}$ is minimal.

Let $\psi:\Delta\to\cO$ be a conformal equivalence with $\psi(0)=z_0\in\cO$,
and let $h=G^+\circ\psi$.  We wish to show that $h$ is a multiple of the
Poisson kernel $P(z,e^{i\kappa})$ for some real $\kappa$.  Let $\lambda$ denote
the measure in the Herglotz representation.  We will show that
the support of $\lambda$ is a single point.   For $0<r<\infty$, let $\cE_r$ be
the unique growth end.  Let $\hat\omega_r=\psi^{-1}(\cO-\bar\cE_r)$, and let
$\omega_r$ be the component of $\hat\omega_r$ containing 0.  Thus
$\Delta\cap\partial\hat\omega_r$ is a collection of Jordan arcs with their
boundaries in
$\partial\Delta$.  Let $\gamma_r$ denote the Jordan arc which separates
$\psi^{-1}(\cE_r)$ from 0.

Let us write $\Delta\cap\partial\omega_r=\gamma_r\cup\bigcup_j\sigma_r^j$. 
Let $\omega_r^j$ denote the component of $\Delta-\sigma_r^j$ which does not
contain 0.  Then $\psi(\hat\omega_r\cap\omega_r^j)$ is a decay end of $\cO$. 
Since there is only one growth end, it follows from Proposition 4.6 that 
$$\lim_{\zeta\in\omega_r^j,\zeta\to\partial\Delta}h(\zeta)=0.$$
Thus $\lambda$ is zero outside the region of $\partial\Delta$ cut out by
$\gamma_r$.  Arguing as in Proposition 4.6, we see that the harmonic length
of $\psi(\gamma_r)$ with respect to $z_0$ inside $\cO$ is  bounded above by
$4\sqrt{|z_0|/r}$.  Transfering this result back to $\Delta$ via
$\psi$, we see that the endpoints of $\gamma_r$ in $\partial\Delta$ are
separated by at most $(2/\pi)\sqrt{|z_0|/r}$. Since the family of curves
$\{\gamma_r:0<r<\infty\}$ is nested, i.e.\ if $r<s<t$, then $\gamma_s$
separates $\gamma_r$ from $\gamma_t$, it follows that they must decrease down to
a single point, which must be the support of
$\lambda$.
\qed

\give Proof of Theorem 4.1.  This is an immediate consequence of Theorems 4.11
and 4.12.

\section 5.  Connectivity of $J$

This Section is devoted to proving 
(Theorem 5.1)  that the presence of either stable or unstable connectivity is
equivalent to the connectivity of
$J$.  

\proclaim Theorem 5.1.  $J$ is connected if and only if either $f$ or $f^{-1}$
is unstably connected.

\proclaim Lemma 5.2.  Let $\nu$ be a hyperbolic measure.  Then for $\nu$ almost
every point $p$ each component of $U^+\cap W^u(p)$ contains $p$ in its closure.

\give Proof.  Let $e(p,r)$ be the number of connected components of $W^u(p)\cap
U^+$ that meet the closed disk of radius $r$ in $W^u(p)$.  Note that if $r<s$
then every component which intersects the disk of radius $r$ also intersects
the disk of radius $s$.  Thus $e(p,r)$ is an increasing function of $r$.  Also,
since every component intersects some disk, we have
$\lim_{r\to\infty}e(p,r)=c(p)$, where $c(p)$ is the number of components of
$W^u(p)\cap U^+$.  By ergodicity, $c(p)$ is constant $\nu$ almost everywhere.

Now $e(fp,\Vert Df_p\Vert)=e(p,1)$ so that
$$e(p,1)=e(fp,\Vert Df_p\Vert)\ge e(fp,1).$$
thus $e(p,1)$ behaves like $c(p,1)$, and arguing as before, we see that for
$\nu$ almost every point $p$ we have $e(p,r)=c(p)$ for every $r$.  In
particular, the set of components that meet the disk of radius $r$ has the same
cardinality as the set of all components.  Thus every component meets the disk
of radius $r$.\qed

\proclaim Lemma 5.3.  Let $\nu$ be a hyperbolic measure.  If $f$ is unstably
connected with respect to $\nu$, then for $\nu$ a.e.\ $p$ the set
$\{G^+<\epsilon\}-\hbox{\rm int}\{G^+=0\}\subset W^u(p)$ is connected.

\give Proof.  The set of points $p$ for which Lemma 5.2 holds has full measure
for $\nu$.  Let $p$ be such a point.  Then for every connected neighborhood
$N$ of $p$ in $W^u(p)$, $N\cup(W^u(p)\cap U^+)$ is connected.  If
$\cO$ is a component of $W^u(p)\cap U^+$, then by Theorem 2.1, $\cO$
satisfies (\dag), so there is a conformal equivalence taking $\cO$ conformally
to the upper half plane and taking $G^+$ to the function $y$.  Thus the
sublevel set $\cO\cap\{G^+<\epsilon\}$ is a strip, which is connected.  Taking
the union over all components, we see that 
$N\cup (W^u(p)\cap \{G^+<\epsilon\})=N\bigcup(\cO\cap\{G^+<\epsilon\})$ is
connected.  Since this holds for all
$N$, the Lemma follows.  \qed

The following general lemma will be used to show that a certain bidisk contains
points of $J$.

\proclaim Lemma 5.4.  Suppose that $(x,y)$ is a local coordinate system such
that $D=\{|x|,|y|<1\}$ has nonempty intersection with both $J^+$ and $J^-$. 
If $G^+>0$ on $\{|x|<1,|y|=1\}$ and $G^->0$ on $\{|x|=1,|y|<1\}$, then it
follows that $\int_D\mu^+\wedge\mu^->0$.

\give Proof. Since
$G^\pm\ge0$  the condition on the supports of $\mu^\pm$ allow us to choose
$\epsilon^\pm>0$ sufficiently small that
$D\cap\{G^+<\epsilon^+\}\subset\{|x|<1,|y|<r\}$ and 
$D\cap\{G^-<\epsilon^-\}\subset\{|x|<r,|y|<1\}$ for some $0<r<1$.
We may assume that $\epsilon^+>0$ is a regular value of $G^+$, so that
$dG^+\ne0$ on $D\cap\{G^+=\epsilon^+\}$.  Thus we may apply Stokes' Theorem:
$$\eqalign{
\int_{D\cap\{G^+<\epsilon^+,G^-<\epsilon^-\}}\mu^+\wedge\mu^-
=&\int_{D\cap\{G^+<\epsilon^+,G^-<\epsilon^-\}}d(d^cG^+\wedge dd^cG^-)\cr
=&\int_{D\cap\{G^-<\epsilon^-\}\cap\partial\{G^+<\epsilon^+\}}d^cG^+\wedge
dd^cG^-.\cr}$$

Since $\epsilon^+$ is a regular value, we may consider
$D\cap\{G^-<\epsilon^-\}\cap\partial\{G^+<\epsilon^+\}$ as a domain inside the
smooth manifold $D\cap\{G^+=\epsilon^+\}$.  It is relatively compact since
$D\cap\{G^+<\epsilon^+,G^-<\epsilon^-\}\subset\{|x|,|y|<r\}$.  Since $G^+$ is
pluriharmonic in a neighborhood of $D\cap\{G^+=\epsilon^+\}$, the closure of
$D\cap\{G^+=\epsilon^+\}$ is not relatively compact in $D$.  Thus $G^-$ is not
constant on this set, so we may choose $\epsilon^-$ to be a regular value for
the restriction of the map $G^-|_{\{G^+=\epsilon^+\}}$.  This means that
$dG^+\wedge dG^-\ne0$ at all points of
$D\cap\{G^+=\epsilon^+,G^-=\epsilon^-\}$, so this is a smooth 2-manifold. 
Since $dd^cG^+=0$ on $\{G^+=\epsilon^+\}$, it follows that $d^cG^+\wedge
dd^cG^-=d(d^cG^-\wedge d^cG^+)$ on
$D\cap\{G^-<\epsilon^-\}\cap\{G^+=\epsilon^+\}$, so we may apply Stokes' Theorem
again to obtain:
$$\eqalign{
\int_{D\cap\{G^-<\epsilon^-\}\cap\partial\{G^+<\epsilon^+\}} d^cG^+\wedge
dd^cG^-
 &=\int_{D\cap\{G^-<\epsilon^-\}\cap\partial\{G^+<\epsilon^+\}}
d(d^cG^-\wedge d^cG^+)\cr
&=\int_{D\cap\partial(\{G^-<\epsilon^-\}\cap\partial\{G^+<\epsilon^+\})}
d^cG^-\wedge d^cG^+.\cr}$$

This last integral is taken over the set $\{G^+=\epsilon^+,G^-=\epsilon^-\}$,
which is an oriented 2-manifold.  The tangent space of a 2-manifold is either a
complex subspace of $\cx2$ (which is equivalent to being invariant under the
complex structure operator $J$) or totally real (which means that $T\oplus
JT=T\cx2$, the generic case).  Since this 2-manifold is compact, it cannot be a
complex submanifold.  Since it is real-analytic, there is an open dense subset
of points where the tangent space it totally real.  The 2-form
$dG^+\wedge dG^-$ annihilates the tangent space.  And since $d^cG^\pm$
is obtained from $dG^\pm$ by applying $J$, it follows that $d^cG^+\wedge d^cG^-$
does not annihilate of the tangent space, so $dG^+\wedge dG^-\wedge
d^cG^+\wedge d^cG^-\ne0$.  It follows that
$dG^+\wedge d^cG^+\wedge dG^-\wedge d^cG^-$ is a nonzero, positive multiple of
the standard volume form on $\cx2$.  From this, we conclude that $d^cG^+\wedge
dG^-\wedge d^cG^-=dG^-\wedge d^cG^-\wedge d^cG^+$ is a nonzero positive multiple
of the volume form on $D\cap\partial\{G^+<\epsilon^+\}$, with the induced
(boundary) orientation; and $d^cG^-\wedge d^cG^+$ is a nonzero, positive
multiple of the volume form on
$D\cap\partial(\{G^-<\epsilon^-\}\cap\partial\{G^+<\epsilon^+\})$ with the
induced orientation.  It follows that the last integral above, and thus
$\int_{D}\mu^+\wedge\mu^-$, is strictly positive. \qed

\give Proof of Theorem 5.1.  We begin by showing that if $f$ is either stably
or unstably connected, then $J$ is connected.  Replacing $f$ by $f^{-1}$, we
may assume that $f$ is unstably connected.  According to Proposition 2.3 of
[BS3], $J^*$ intersects every connected component of $J$.  Periodic saddle points
are dense in $J^*$, by [BLS], so it suffices to show that any two saddle points
$p$ and $q$ can be connected by a path lying in an arbitrarily small
neighborhood $U$ of $J$.  Again by [BLS], any two saddle points are heteroclinic,
so $W^u(p)\cap W^s(q)$ is nonempty and in fact contains a transverse
intersection.  Now $W^u(p)\cap W^s(q)$ contains points arbitrarily close to $q$,
so it suffices to show that any transverse intersection point $r\in W^u(p)\cap
W^s(p)$ can  be connected to $p$ by a path lying inside  $U$.  

It is evident that $r\in W^u(p)\cap\{G^+=0\}$, and we will show that
$r\notin\hbox{\rm int}\{G^+=0\}$, where the interior is taken relative to
$W^u(p)$.  Let us suppose, to the contrary, that there is a disk in
$\{G^+=0\}\cap W^u(p)$ containing $r$.  Since the iterates of $f^n$, $n\ge0$
remain bounded, the derivative $Df^n$ tangential to $W^u(p)$ at $r$ remains
bounded.  But this contradicts the smooth Lambda Lemma.  Thus $r$ is not in the
interior, and $r$ must belong to the closure of $\{0<G^+<\epsilon\}\cap W^u(p)$
inside $W^u(p)$, and so it will follow from Lemma 5.3 that there is a path in
$\{G^+<\epsilon\}\cap W^u(p)$ connecting $r$ to $p$.  Since we may choose
$\epsilon>0$ sufficiently small that $\{G^+<\epsilon\}\cap W^u(p)\subset U$, it
follows that $r$ and $p$ are in the same connected component of $J$.  This
completes the proof that $J$ is connected.

Now we show that if neither $f$ nor $f^{-1}$ is unstably connected, then $J$ is
not connected.  Let $p$ be a periodic saddle point.  Replacing $f$ by $f^n$ for
an appropriate $n$ lets us assume that $p$ is a fixed point.  We can choose a
coordinate system $\{|x|<1,|y|<1\}$ in a neighborhood $B$ of $p$ so that $p$
corresponds to the point (0,0), the set $\{|x|<1,y=0\}$ is a local unstable
manifold for $p$, and the set $\{x=0,|y|<1\}$ is a stable manifold for $p$. 
Furthermore, by taking $B$ small, we may assume that the restriction of $f$
to $B$ is approximately linear, and thus it is uniformly expanding in the
$x$-direction and uniformly contracting in the $y$-direction.

Since $f$ is not unstably connected, $W^u(p)\cap K^+$ contains a compact
component.  Now $f^{-1}$ decreases distance in $W^u(p)$, so by applying
$f^{-m}$ with $m$ sufficiently large we may assume that the set
$\{|x|<1,y=0\}\cap K^+$ has a compact component.  We may further assume that
this component does not contain $p=(0,0)$.  Let $\gamma $ be a curve in
$\{|x|<1,y=0\}\cap U^+$ which encloses this compact component but does not
enclose $p$. Let $D_0$ be the region of $\{|x|<1,y=0\}$ enclosed by $\gamma$,
and let $E^+$ denote the portion of $W^u(p)\cap K^+$ enclosed by $\gamma$. 
By the uniform expansion/contraction of $f|_B$ (or equivalently, by the
Lambda Lemma) it follows that for $n$ sufficiently large and $|x|<{1\over2}$, the
disk
$$M^+_n(x):=f^{-n}(\{x\}\times\{|y|<\epsilon\})\cap\bigcap_{j=0}^{n-1}f^{-j}B$$ 
is vertical, in the sense that the projection of $M^+_n(x)$ to
$\{x=0,|y|<\epsilon\}$ is a homeomorphism.  It follows that
$$\Gamma^+_n:=f^{-n}(\gamma\times\{|y|<\epsilon\})
\cap\bigcap_{j=0}^{n-1}f^{-j}B=\bigcup_{x\in\gamma}M^+_n(x)$$
is a hypersurface in $B$ made up of vertical disks.  Thus $\Gamma^+_n$ divides
$B$ into two components, and we let $B^+$ denote the component of
$B-\Gamma^+_n$ which contains $f^{-n}D_0$.  We note, also, that $G^+>0$ on
$\Gamma^+_n$.

Similarly, since $f^{-1}$ is not unstably connected, $W^s(p)\cap K^-$ contains
a compact component.  Arguing as above, we have a hypersurface $\Gamma^-_n$ of
unstable disks $M^-_n(x)$, and $G^->0$ on $\Gamma^-_n$.  Furthere, there is a
component $B^-$ of $B-\Gamma^-_n$ which contains the compact component of
$W^s(p)\cap K^-$.  Let $B_0=B^+\cap B^-$.  We consider the family of vertical
disks
$\{|x|<{1\over2}\}\ni x\mapsto M^+_n(x)$.  We know that $G^->0$ on $\Gamma^-_n$,
and so $dd^cG^-$ vanishes there.  Thus  $\{|x|<{1\over2}\}\ni x\mapsto
\int_{B^-\cap M^+_n(x)} dd^cG^-$ is constant.  We know that
$B^-\cap M^+_n(0)=B^-\cap W^s(p)$ contains a compact component, so it follows 
$dd^cG^-$ puts positive mass on $B^-\cap M^+_n(0)$.  Thus $dd^cG^-$ puts
positive mass on $B^-\cap M^+_n(x)$ for $x\in D_0\subset\{|x|<{1\over2}\}$.  This
implies that
$B_0$ intersects $J^-$.  Similarly, $B_0$ intersects $J^+$.  Thus we may apply 
Lemma 5.4 to conclude that $\int_{B_0}\mu^+\wedge\mu^-\ne0$, from which it
follows that $J\cap B_0\ne\emptyset$.  Since $\partial
B_0\subset\Gamma^+_n\cup\Gamma^-_n$ is disjoint from $J$, and since $B_0$ does
not contain $p$, and thus all of $J$, it follows that
$J$ is disconnected. \qed

\section 6. Unstable connectivity and extension of $\varphi^+$

In this section we find several characterizations of the condition of unstable
connectivity which are summarized in Theorem 6.3. These relate to the existence of 
extensions of
$\varphi^+$ and topological properties of $J^-_+$.

Recall that $\varphi^+$ is defined and holomorphic on $V^+$ and satisfies the
functional equation $$\varphi^+(f^n(p))=(\varphi^+(p))^{d^n}\eqno{(6.1)}.$$ When $f$
is unstably connected then according to Theorem 2.1 the function $\varphi^+$ has a
continuous extention to  $J^-_+$ that satisfies equation (6.1).
Let $\cal G^+$ be the foliation of $U^+$ defined by the holomorphic one
form $\partial G^+$. When restricted to $V^+$ the leaves of $\cG^+$ are just the sets
on which $\varphi^+$ is constant.

\proclaim Lemma 6.1. If $f$ is unstably connected then each path component of $J^-_+$
is simply connected.

\give Proof.  Let $\gamma$ be a loop in $J^-_+$. The image $\varphi^+(\gamma)$ is a loop
in $\C-\bar\Delta$. We begin by showing that the image loop is contractible. Let
$[\varphi^+(\gamma)]\in\Z$ denote the degree of the image loop. To show that the image
loop is contractible we show that the degree is zero. Let $\gamma'$ denote
$f^{-n}(\gamma)$. Now the functional equation for $\varphi^+$ gives:
$$[\varphi^+(\gamma)]=[\varphi^+(f^n(\gamma')]=[(\varphi^+)^{d^n}(\gamma')]=
d^n[\varphi^+(\gamma')].$$
This implies that $[\varphi^+(\gamma)]/d^n$ is an integer for any $n$. Thus
$[\varphi^+(\gamma)]=0$. 

Lemma 2.4 tells us that the map $\varphi^+: J^-_+\to\C-\bar\Delta$ is a locally trivial
fibration. Thus it has the homotopy lifting property. Since the image of $\gamma$ is
contractible the loop $\gamma$ is homotopic to a loop in a fiber of the map $\varphi^+$,
that is to say a set $\varphi^+=const$. We complete the proof by showing that each
component of a fiber is simply connected. Now a fiber of the map is contained in a leaf
$L$ of the foliation
$\cG^+$. In fact it is contained in the intersection of $L$ with $J^-$. The set $J^-$ is the
zero set of the function $G^-$. The leaf
$L$ is conformally equivalent to $\C$ and the restriction of $G^-$ to $L$ is a subharmonic
function. The maximum principle implies that each component of the zero set of $G^-$ in $L$
 is simply connected.  \qed

\proclaim Lemma 6.2.  If $f$ is unstably connected, then $\varphi^+$ has an
analytic continuation to a neighborhood of $J^-_+$.

\give Remark. It follows from Hubbard and Oberste-Vorth [HO] that $\varphi^+$
cannot be extended to $U^+$.  Any holomorphic extension of $\varphi^+$ is locally 
constant on the leaves of $\cal G^+$.  Thus if an  extension of $\varphi^+$ to a set
$U'\supset V^+$ exists, each leaf of $\cG^+|U'$ can intersect only one disk of $\cG^+|V^+$,
since $\varphi^+$ takes distinct values on distinct disks.  There can be no
extension to
$U^+\supset V^+$, since as shown in [HO] each leaf of $\cG^+$ intersects $V^+$ in infinitely
many disks.

\give Proof.    For $p\in\cal G^+$ let $L_p$ denote the leaf of the $\cal G^+$ foliation
that passes through $p$. For $p$ and $q$ in the same leaf let $d_L(p,q)$ denote the
distance measured with respect to the induced Riemannian metric in the leaf. For $p\in
U^+$ let
$\nu_p$ be the set consisting of ``nearest points in $J^-_+$,''  i.e., those points
$q$ in $L_p\cap J^-_+$  which minimize the function
$d_L(p,q)$ among all points in $L_p\cap J^-_+$. Let $N$ consist of those points $p$ for
which the function
$\varphi^+$ is constant on
$\nu_p$. For $p\in N$ define the function $\tilde\varphi(p)$ to be the common value of
$\varphi^+$ on the elements of $\nu_p$. We will show that any $p\in J^-_+$ has a
neighborhood in $N$ on which the function
$\tilde\varphi$ is holomorphic. Choose  local coordinates $u$ and $v$ near $p$ so
that the set $B=\{(u,v): |u|\le1,\ |v|\le1\}$ is a neighborhood of $p$, and the sets
$u=const$ are contained in leaves of $\cal G^+$. We may assume that the set $v=0$ is
the local leaf of the $J^-$ lamination containing $p$. Choose an
$n$ sufficiently large so that
$f^n(p)\in V^+$. We may assume that $B$ is chosen small enough so that $f^n(B)\subset V^+$.
For $(u,v)\in J^-_+\cap B$ define the following function
$\alpha(u,v)=\varphi^+(u,v)/\varphi^+(u,0)$. 
Let $W$ denote the set where $\alpha=1$. Since $\alpha$ is continuous $W$ is a closed
set. We claim that $W$ is an open set of $J^-_+\cap B$.

Since $(u,v)$ and $(0,v)$ have the same second coordinate they lie in a disk inside a leaf
of the $\cG^+$ foliation. Since $f^n(B)\subset V^+$ the points
$f^n(u,v)$ and
$f^n(0,v)$ are on the same leaf of the $\cG^+$ foliation of $V^+$ and we have
$\varphi^+(f^n(u,v))=\varphi^+(f^n(0,v))$. The functional equation (6.1) gives
$$(\varphi^+)^{d^n}(u,v)=(\varphi^+)^{d^n}(u,0).$$ So $\alpha(u,v)^{d^n}=1$ and the values
of $\alpha$ are $d^n$-th roots of unity. The function $\alpha$ is continuous and takes on a
finite set of values so the set where $\alpha=1$ is open.   The set
$\{(v,0)\}$ is in $W$. Since $W$ is open we can choose $\epsilon$
sufficiently small so that
$|u|<\epsilon$ implies that the $d_L$ distance from $(u,v)$ to $(u,0)$ is smaller
than the $d_L$ distance from $(u,v)$ to any point $(u,v')$ not in $W$ or any point
$(u,v')$ on the boundary of $B$. Since the nearest neighbors of $(u,v)$ are in $W$ we
have that $\varphi^+(u,v')=\varphi^+(u,0)$ for all nearest neighbors of $(u,v)$. Thus for
$|v|<\epsilon$
$(u,v)\in N$ and $\tilde\varphi(u,v)=\varphi^+(u,0)$.  \qed

\proclaim Theorem 6.3. The following are equivalent:
\item{(1)} $f$ is unstably connected.
\item{(2)} $\varphi^+$ extends to a continuous function on $J^-_+$ which satisfies the
equation $\varphi^+f^n=(\varphi^+)^{d^n}$.
\item{(3)} $\varphi^+$ extends to a continuous function on $J^-_+$ which is
holomorphic on leaves.
\item{(4)} The cohomology class represented by the form $\eta=(1/2\pi)d^cG^+$ is an
integral class on each leaf of the lamination $\cM^-$ of $J^-_+$.
\item{(5)} Each path component of $J^-_+$ is simply connnected.
\item{(6)} $H_1(J^-_+;\R)=0$.
\item{(7)} $\varphi^+$ extends holomorphically to a neighborhood of $J^-_+$.

\give Proof. The strategy of proof is to show that (1)$\Rightarrow$(2),
(2)$\Rightarrow$(3), (3)$\Rightarrow$(4) and (4)$\Rightarrow$(1).
We then show that (1)$\Rightarrow$(5),
(5)$\Rightarrow$(6) and (6)$\Rightarrow$(4). We conclude by showing that
(1)$\Rightarrow$(7) and (7)$\Rightarrow$(3).

(1)$\Rightarrow$(2). This follows from Theorem 2.1. 

(2)$\Rightarrow$(3) The function $\varphi^+$ is defined and holomorphic on $V^+$. Let
$p\in J^-_+$. For some $n$, $f^n(p)\in V^+$. Now the function $\varphi^+f^n$ is
holomorphic when restricted the leaf containing $p$. The extension is locally a
continuous $d^n$-th root of a holomorphic function. Hence the extension is holomorphic on
each leaf.

(3)$\Rightarrow$(4). Since $\varphi^+$ is holomorphic on leaves the function
$\log|\varphi^+|$ is harmonic on leaves. In the set $V^+$ we have $\log|\varphi^+|=G^+$. 
Since both sides of the equation are analytic functions the equation holds on the entire
leaf. Now 
$$\eta=(1/2\pi) d^cG^+=(1/2\pi) d^c\log|\varphi^+|=(\varphi^+)^*((1/2\pi) d^c\log|z|).$$
The form $d^c\log|z|$ measures the change in argument so on any closed loop its value is
in the set $2\pi\Z$. Thus $\eta$ is the pullback of a form which represents an integral
class in
$H^1(\C-\bar\Delta)$ so $\eta$ itself is an integral class.

(4)$\Rightarrow$(1).
Let $p$ be a saddle point. Assume $\eta$ represents an integral class but that $f$ is 
unstably disconnected. Let $\cO$ be a component of $W^u(p)-K^+$.  By Theorem 4.11 there
are only finitely many components, so we may assume that $\cO$ is periodic. 
Passing to a power of $f$, we may assume that $\cO$ is fixed.  Since $f$ is unstably
disconnected there is a nonempty compact component $E$ of $W^u(p)-\cO$. Let 
$\gamma\subset\cO$ be a simple closed curve that surrounds $E$. Let $D$ be the topological
disk surrounded by $\gamma$.  Define
$\delta_\gamma:=\int_\gamma d^cG^+$.  Given an $n$ let $\gamma'=f^{-n}(\gamma)$. 
Our hypothesis implies that $\delta_\gamma$ and $\delta_\gamma'$ are
integers. Now the functional equation for $G^+$ gives:
$$\delta_{\gamma}=\int_\gamma
d^cG^+=\int_{\gamma'}d^cG^+\circ f^n=d^n \int_{\gamma'}d^cG^+ = d^n\delta_\gamma'$$ Since
$\delta_\gamma/d^n$ is an integer for any $n$, $\delta_\gamma=0$.

  Since $G^+$ is subharmonic we have $\int_E dd^cG^+\ge0$. If $\int_E dd^cG^+\ge0=0$ then
$G^+$ would be harmonic on the region enclosed by $\gamma$. Since $G^+$ is positive on
$\gamma$ and zero on $E$ this would violate the minimum principle. We conclude that
$\int_E dd^cG^+>0$. But $\int_E dd^cG^+=\delta_\gamma=0$. This contradiction completes
the proof. \qed

(1)$\Rightarrow$(5). This is Theorem 6.1.

(5)$\Rightarrow$(6). This is clear.

(6)$\Rightarrow$(4). If  $H_1(J^-_+;\R)=0$ then the simplicial cohomology group
$H^1(J^-_+;\R)$ is zero. Since $\eta$ represents an element of this group $\eta=0$. In
particular $\eta$ is integral.

(1)$\Rightarrow$(7). This is Theorem 6.2.

(7)$\Rightarrow$(3). This is clear. \qed

\section 7.  Critical Points and Harmonic Measure 

We give first a dichotomy of possible behaviors: either $f$ is unstably
connected, or $f$ has a  strong unstable disconnectedness property with respect to
harmonic measure.  Then we will show that $f$ is unstably connected exactly when it has no
critical points with respect to $\mu$. 

Let $\cR$ denote the set of Pesin regular points for the map $f$.  For
$p\in\cR$ we let $K^{+,u}(p)$ denote the connected component of $W^u(p)\cap K^+$
which contains $p$.  We let $\cT^u$ denote the set of points of $\cR$ for which
the corresponding component of the unstable slice of $K^+$ is trivial, i.e.
$\cT^u=\{p\in\cR:\{p\}=K^{+,u}(p)\}$.

\proclaim Theorem 7.1.  The following dichotomy holds.  Either
\item{1.} $f$ is unstably connected, or
\item{2.} $\cT^u$ has full measure for $\mu$. Equivalently, for $\mu$ a.e.\ $p$,
$\cT^u\cap W^u(p)$ has full measure for the induced measure $\mu^+|_{W^u(p)}$. 
\vskip0pt\noindent In case (2), it follows that for $\mu$ a.e.\ point $p$, the critical
points of $G^+|_{W^u(p)\cap U^+}$ are dense in the boundary of $W^u(p)\cap
K^+$, with the boundary being taken inside $W^u(p)$.

\give Remark.  This dichotomy says that, with respect to harmonic measure,
either $f$ is unstably connected, or unstably totally disconnected.  We recall
that the induced measure $\mu^+|_{W^u(p)}$ is defines the measure class of
harmonic measure.  We conclude that when case (2) of Theorem 6.2 holds, and
there is a nontrivial compact component $E$, then $E$ has zero harmonic
measure.  In this case the point components of
$W^u(p)\cap K^+$ form a halo which surrounds $E$ closely enough that $E$ can
have no harmonic measure.

\give Proof.  For $p\in\cR$, $W^u(p)\approx\cx{}$ has an affine structure, and
we let $\Vert\cdot\Vert_G$ denote the norm defined in \S2.  Let us define
$R(p)$ to be the radius with respect to $\Vert\cdot\Vert_G$ of the
smallest closed disk centered at $p\in W^u(p)$ which contains $K^+(p)$.  It
follows from the ergodicity of $\mu$ that either: (1) $R(p)=\infty$, $\mu$
a.e., (2) $R(p)=0$, $\mu$ a.e., or (3) $0<R(p)<\infty$, $\mu$ a.e.  The
possibility (2) corresponds to the case (2) in the
dichotomy above.  

We show that case (1) here corresponds to case (1) above.  Let us find countably
many unstable boxes $B_j^u$ whose union has full measure.  If $f$ is not unstably
connected, we may choose one of the unstable boxes $B_j^u=\{\Gamma(t): t\in
T_j\}$ so that one of the leaves
$\Gamma(t_0)$ intersects $K^+$ in a compact component.  Thus there is a
simple, closed curve $\gamma\subset\Gamma(t_0)$ which encircles a nonempty
portion of $K^+\cap\Gamma(t_0)$.  It follows that $\min_\gamma G^+>0$, and
by continuity $\min_{\gamma_t}G^+>0$ for $\gamma_t\subset\Gamma(t_0)$ near
$\gamma$.  Further, it follows that $t\mapsto\int_{\gamma_t}d^cG^+$ is locally
constant for $t$ near $t_0$.  Thus $\gamma_t$ cuts off a compact portion $E_t$
of $K^+\cap \Gamma(t)$ for $t$ near $t_0$.  The set of $p$ such that 
$W^u(p)\cap K^+$ has a compact component contains $\bigcup_{|t-t_0|<\epsilon}
E_t$ which has positive $\mu$ measure by the local product structure, since the
measure of $E_t$ is  $\mu^+|_{\Gamma(t)}(E_t)=\int_{\gamma_t}d^cG^+>0$, and
$\{t\in T:|t-t_0|<\epsilon\}$ has positive transversal measure.  Thus
$f$ is not unstably connected.

Thus to prove the Theorem, we must show that (3) cannot occur.  For
$0<a<b<\infty$ we define $S=\{p\in\cR:a<R(p)<b\}$.  In case (3), we may choose
$a$ and $b$ so that $S$ has positive $\mu$ measure.  Now by (4.1) $f^n$ is
linear with respect to the affine structures of $W^u(p)$ and $W^u(f^np)$.  Thus
it follows that
$$R(f^np)=\Vert Df^n|_{E^up}\Vert_G R(p).$$
By Poincar\'e recurrence, for almost every $q\in S$, there is a
sequence $n_j\to\infty$ such that $f^{n_j}q\in S$.  In this case we have
$b/a>\Vert Df^{n_j}|_{E^up}\Vert_G$.  But this contradicts Lemma 4.4 if
$S$ has positive measure.  Thus case (3) cannot occur.

The equivalent statement in (2) follows by the local product structure, which
says that $\mu$ is given locally as the product of the slice measures of $\mu^+$
and $\mu^-$.  Thus if $\cT^u$ has full measure for $\mu$, then it has full
measure for almost every slice measure $\mu^+|_{W^u(p)}$.  

For the statement concerning critical points, 
we note that if $q\in\cT^u\cap W^u(p)$, then $G^+|_{W^u(p)\cap U^+}$ must have
critical points arbitrarily close to $q$.  Otherwise, by Lemma 6.2, $q$ is an
isolated point of $W^u(p)\cap K^+$, which is impossible since it would mean
that $q$ is an isolated zero of $G^+|_{W^u(p)}$.
\qed

\proclaim Lemma 7.2.  If $E$ is a compact component of
$W^u(p)\cap K^+$, and if $G^+|_{W^u(p)}$ has no critical points in a
neighborhood of $E$, then $E$ is an isolated component of $W^u(p)\cap K^+$.

\give Proof.  Let $V$ denote a relatively compact neighborhood of $E$ inside
$W^u(p)$ such that $\partial V\cap K^+=\emptyset$.  If we
set $\delta_0=\min_{\partial V}G^+$, then the sublevel sets
$S_\delta:=\{G^+<\delta\}\cap V$ are bounded if
$0<\delta<\delta_0$, and $\partial S_\delta\subset\{G^+=\delta\}$.  Since $G^+$
has no critical points, the set $\{G^+=\delta\}$ is smooth, and each component of
$\partial S_\delta$ is homeomorphic to a 1-sphere.  It follows by Morse theory
that $E$ can be the only component of $K^+$ inside its component of $S_\delta$. 
Thus $E$ is isolated.  \qed

\proclaim Theorem 7.3.  The following are equivalent:
\item{(1)} $\lambda^+(\mu)=\log d$
\item{(2)}  For $\mu$ a.e.\ $p$, $G^+|_{W^u(p)-K^+}$ has no critical points.
\item{(3)} For a set of $p$ of positive
$\mu$ measure there is a component $\cO_p$ of $W^u(p)-K^+$ such that
$G^+|_{\cO_p}$ has no critical points. 
\item{(4)} $f$ is unstably connected.

\give Proof.  $(1)\Leftrightarrow(2)$ follows from Corollary 6.7 of
[BS5].  $(2)\Rightarrow(3)$ is obvious. $(4)\Rightarrow(2)$ follows from Corollary 2.19.  It
remains to show $(3)\Rightarrow(4)$.   By Theorem 6.1, there are two possibilities: if $f$
is not unstably connected, then for almost every $p$, the unstable manifold $W^u(p)$ has the
property that $\mu^+|_{W^u(p)}$ is carried by the set $\cT^u$.  But if $q\in\cT^u$, then the
$\{q\}$ is the component of
$W^u(p)\cap K^+$ containing $q$, which must be isolated by Lemma 6.2, since
$G^+|_{W^u(p)}$ has no critical points.  Since $G^+$ is continuous, however,
$\mu^+|_{W^u(p)}$ can put no mass on an isolated point.  Thus the slice cannot
put any mass at all on
$\cT^u$, so $f$ must be unstably connected. \qed

\proclaim Corollary 7.4.  If $f$ is dissipative then  $f$ is unstably connected.  If $f$
preserves volume, then it is stably connected if and only if it is unstably connected.

\give Proof.  By Theorem 7.3, if $f$ is unstably connected, then $\Lambda=\log
d$.  If $a\in\cx{}$ denotes the (constant) complex jacobian determinant of $f$,
then by Proposition 7.7 of [BS5], $\Lambda=\log d$ implies that $|a|\le 1$. 
Similarly, if $|a|=1$, then $\Lambda(f)=\log d$ if and only if
$\Lambda(f^{-1})=\log d$.  \qed

\bigskip
\centerline{\bf References}
\bigskip

\item{[BLS]} E.\ Bedford, M.\ Lyubich, and J.\ Smillie, Polynomial diffeomorphims
of ${\bf C}^2$. IV: The measure of maximal entropy and laminar currents.  Invent.
math. 112, 77--125 (1993).

\item{[BLS2]} E.\ Bedford, M.\ Lyubich, and J.\ Smillie, Distribution of periodic points of
polynomial diffeomorphisms of $\C^2$.  Invent. math. 114, 277--288 (1993).

\item{[BS3]} E. Bedford and J. Smillie, Polynomial diffeomorphisms of ${\bf C}^2$.
III: Ergodicity, exponents and entropy of the equilibrium measure. Math. Ann.
294. 395--420 (1992).

\item{[BS5]} E. Bedford and J. Smillie, Polynomial diffeomorphims
of ${\bf C}^2$. V: Critical points and Lyapunov exponents. J. of Geometric
Analysis, to appear.

\item{[BS7]} E. Bedford and J. Smillie, Polynomial diffeomorphims
of ${\bf C}^2$. VII: Hyperbolicity and external rays.

\item{[C]}  A. Candel,  Uniformization theorem for surface laminations, Ann. Scient
\'Ec. Norm Sup., 4 s\'erie t. 26 (1993), 489-516.

\item{[DH]}  A. Douady and J. Hubbard, It\'eration des polyn\^omes quadratiques complexes, C.
R. Acad. Sc., 294 s\'erie I (1982), 123-126.

\item{[FS]}  J.-E. Fornaess and N. Sibony, Complex dynamics in higher
dimensions, in {\sl Complex Potential Theory}, 131--186, P. Gauthier (ed.) 1994.

\item{[F]}  W.H.J. Fuchs, {\sl Topics in the Theory of Functions of One Complex
Variable}, van Nostrand, 1967.

\item{[H]} J.H. Hubbard, Local connectivity of Julia sets and bifurcation loci:
Three theorems of J.-C. Yoccoz, in {\sl Topological Methods in Modern
Mathematics}, 467--511, L. Goldberg and A. Phillips, (ed.) 1993.

\item{[HO]} J.H. Hubbard and R. Oberste-Vorth, H\'enon mappings in the complex
domain I: The global topology of dynamical space, Inst.\ Hautes \'Etudes Sci.\
Publ.\ Math.\ 79, 5--46 (1994).

\item{[KM]}  A. Katok and L. Mendoza, Dynamical systems with nonuniformly
hyperbolic behavior, supplement to {\sl Introduction to the Modern Theory of
Dynamical Systems}, A. Katok and B. Hasselblatt, Cambridge University Press,
1995.

\item{[LS]} F. Ledrappier and J.-M. Strelcyn, A proof of the estimation from
below in Pesin's entropy formula, Ergod.\ Th.\ \&\ Dynam.\ Sys.\ 2, 203--219
(1982).
 
\item{[MSS]}  R. Man\'e, P. Sad, and D. Sullivan,  On the dynamics of
rational maps,   Ann.\ scient.\ \'Ec.\ Norm.\ Sup. 16, 193--217 (1983).
 
\item{[P]}  M.\ Pollicott,  {\sl Lectures on ergodic theory and Pesin theory
on compact manifolds}, London Mathematical Society Lecture Note Series 180,
Cambridge U. Press, 1993.

\item{[S]} E. Spanier, {\sl Algebraic Topology}, McGraw-Hill, 1966.

\item{[Su]} D. Sullivan, Linking the universalities of Milnor-Thurston Feigenbaum and
Ahlfors-Bers, in {\it Topological methods in modern mathematics}, 543-564,  L. Goldberg
and A. Phillips, (ed.) 1993.

\item{[Wa]} P. Walters, {\sl An Introduction to Ergodic Theory}, Springer-Verlag,
1982.

\item{[W]} H. Wu, Complex stable manifolds of holomorphic diffeomorphisms,
Indiana U. Math.\ J. 42, 1349--1358 (1993).

\end